\documentclass[11pt,a4paper,reqno]{amsart}
\usepackage{epsfig,amssymb,latexsym,amscd}

\usepackage[all,2cell]{xy}
\xyoption{v2}
\UseAllTwocells
\theoremstyle{plain}

\newtheorem{theo}{Theorem}
\newtheorem{coro}{Corollary}
\newtheorem{lema}{Lemma}

\theoremstyle{remark}

\theoremstyle{definition}

\newtheorem{defi}{Definition}
\newtheorem{obse}{Observation}

\newcommand{\compo}{{\scriptstyle \circ}}

\newcommand{\bicomod}{{}^C\!\mathcal M\hspace*{.2mm}^C}
\newcommand{\hopfmod}{{}^C\!\mathcal M^C_C}

\newcommand{\rcomod}{\mathcal M^C}
\newcommand{\bicomodH}{{}^H\!\mathcal M\hspace*{.2mm}^H}
\newcommand{\hopfmodH}{{}^H\!\mathcal M^H_H}
\newcommand{\lcomodH}{{}^H\!\mathcal M}
\newcommand{\rcomodH}{\mathcal M^H}

\newcommand{\id}{\ensuremath{\mathrm{id}}}
\newcommand{\M}{\ensuremath{\mathcal{M}}}

\begin{document}
\title[Monoidal categories of comodules]{Monoidal categories of comodules for coquasi Hopf algebras and
  Radford's formula}
\author{Walter Ferrer Santos}
\address{Facultad de Ciencias\\Universidad de la Rep\'ublica\\
Igu\'a 4225\\11400 Montevideo\\Uruguay\\}
\thanks{The first author would like to thank, Csic-UDELAR, Conicyt-MEC}
\email{wrferrer@cmat.edu.uy}
\author{Ignacio Lopez Franco}
\address{Centre for Mathematical Sciences\\Wilberforce Road\\Cambridge CB3 0WB\\UK}
\thanks{The second author acknowledges the support of a Internal Graduate
Studentship form Trinity College, Cambridge}
\email{I.Lopez-Franco@dpmms.cam.ac.uk}

\begin{abstract}
We study the basic monoidal properties of the category of Hopf modules
for a coquasi Hopf algebra. In particular we discuss 
the so called fundamental
theorem that establishes a monoidal equivalence between the category
of comodules and the category of Hopf modules. 
We present a categorical proof of Radford's
$S^4$ formula for the case of a finite dimensional coquasi Hopf
algebra, by establishing a monoidal isomorphism between certain double
dual functors.
\end{abstract}
\maketitle

\begin{center}
{\em Dedicated to I. Shestakov on the occasion of his 60th birthday}
\end{center}
\section{Introduction}
\label{section:intro}
The main purpose of this paper is to study the monoidal category of
Hopf modules for a coquasi Hopf algebra. As a consequence we obtain 
a proof of Radford's
$S^4$ formula valid for finite dimensional coquasi Hopf algebras. 
Inspired in \cite{kn:eno} we show that this formula is
intimately related to the existence of certain natural transformation 
relating the left and the right double dual functors 
for the category of right $H$--comodules. This natural transformation
comes from the application of the structure theorem for Hopf modules,
to ${}^{*}H$ viewed as a right $H$--Hopf module. 

Coquasi Hopf algebras are the dual notion of the quasi Hopf algebras
defined in \cite{kn:dri}. The main difference with Hopf algebras is
that for coquasi Hopf algebras the role of the multiplicative and
comultiplicative structures is not longer interchangeable. 
In a coquasi Hopf algebra the
multiplicative structure is no longer one dimensional, but two
dimensional; this is expressed in the fact that the multiplication is
not longer associative but only up to isomorphism, provided by a
functional $\phi$. The antipode is also defined as a two dimensional
structure, the extra dimension provided by two functionals
$\alpha,\beta$. See below.

The category of Hopf modules in the context of (co)quasi Hopf algebras
has been considered by
different authors and it was initially studied
in \cite{kn:hn,kn:Schauenburg}.

In the case of Hopf algebras, Radford's formula for $S^4$  
was first proved in full
generality in \cite{kn:r}, with predecessors in \cite{kn:l} and 
\cite{kn:s}. A more recent proof, appears in \cite{kn:schneider}. 
There are many generalizations of the formula
from the case of Hopf algebras to other situations, {\em e.g.}: braided
Hopf algebras,  bF algebras --braided and classical--, quasi Hopf
algebras, weak Hopf algebras, and Hopf algebras over rings. 
The following is a partial list of references for some proofs of these
generalizations:  \cite{kn:bkl}, \cite{kn:dt}, \cite{kn:fs}, \cite{kn:hn},
\cite{kn:ks}, \cite{kn:k} and \cite{kn:n}. 
Closer to the spirit of our paper an analogue of
Radford's formula for finite tensor categories appears in \cite{kn:eno}.

In this Introduction, after a general description of the paper, we
briefly recall the definition and the first basic properties 
of a coquasi bialgebra and of a coquasi Hopf algebra. We define the
notion of monoidal morphism that is the adequate notion of morphism between
coquasi bialgebras. Later in this
Introduction we establish the basic notations that will be used along
the paper.  

In Section \ref{s:catsofbicomod}, that is the technical core of the
paper, we present the basic properties of the monoidal categories used
later. We work with the categories of comodules and Hopf modules 
for coquasi bialgebras and recall the properties of 
the monoidal structures induced by the cotensor product over $H$ and
by the tensor product over $\Bbbk$.
We look at some basic monoidal functors
associated to monoidal morphisms given by corestriction of scalars
and
its adjoints given by coinduction. We also consider 
other monoidal functors --{\em e.g.}, the left adjoint comodule and
right adjoint comodule functors--
that will be used later. 

In Section \ref{s:duality} we recall that in the case of coquasi Hopf
algebras with invertible antipode the tensor 
categories of finite dimensional comodules --or bicomodules-- are
rigid, {\em i.e.}, each object has a left and a right dual. 
We use this rigidity in order to describe
explicitly --for
finite dimensional Hopf algebras-- the monoidal structure of the
antipode. 
 
In Section \ref{s:fundtheo} we present a proof of the version of the
fundamental theorem on Hopf modules for coquasi Hopf algebras 
that we need later in the
paper. By applying some general results on Hopf modules over autonomous
pseudomonoids to our context we prove that the 
free right Hopf module functor is a monoidal
equivalence from the category of comodules into the category of Hopf
modules. 

In Section \ref{s:frobiso}, we apply the fundamental theorem on Hopf modules to
${}^*\!H$ and obtain the Frobenius isomorphism that is a morphism in
the category of Hopf modules between $H$ and ${}^*H$. 
Along the way we identify the one dimensional object of cointegrals in
this context.
 
In Section \ref{s:formula}, using the categorical machinery
constructed above and in the same vein than in \cite{kn:eno}, we prove the
existence of a natural monoidal isomorphism between the double duals on the
left and on the right of a finite
dimensional left $H$--module. This isomorphism 
will yield Radford's formula.

It is not obvious {\em a priori}\/ that the formula obtained for finite
dimensional coquasi Hopf algebras is related to the classical Radford's
formula for Hopf algebras.  
Thus, in Section \ref{s:casehopf}, we apply the previously developed
techniques to the case that $H$ is a classical Hopf
algebra in order to deduce the original Radford's formula for $S^4$
--see \cite{kn:r}--. 

In Section \ref{s:back} we present in an Appendix 
the categorical background needed to 
prove some of the basic monoidal properties of the cotensor
product. We recall a few basic definitions and results about density of
functors and completions of categories under certain classes of
colimits.

\subsection{Basic definitions}
Next we summarize the basic definitions that we need.

Recall that the category of coalgebras and morphisms of coalgebras has
a monoidal structure such that the forgetful functor into the category
of vector spaces is monoidal. In other words, the tensor product over
$\Bbbk$ of two coalgebras is a coalgebra and $\Bbbk$ is a coalgebra, in
a canonical way. This is a consequence of the fact of that  category
of vector spaces is braided, and in fact symmetric, with the usual
switch $\mathrm{sw}:V\otimes W\to W\otimes V$.

Assume $(C,\Delta, \varepsilon)$ is a coalgebra. The maps
$\Delta:C\to C\otimes C$ and $\varepsilon:C\to\Bbbk$ are the comultiplication
and the counit respectively. We will use Sweedler's notation
as introduced in \cite{kn:sweedler}, and write $\Delta(c)=\sum c_1\otimes c_2$.
We use the notation $\Delta^2$ for the morphism $\Delta^2(c)= \sum c_1
\otimes c_2 \otimes c_3$. Moreover, the convolution product will be
denoted by the symbol $\star$.  
\begin{defi}
A {\em coquasi bialgebra}\/ structure on the coalgebra
$(C,\Delta,\varepsilon)$ is a triple
$(p,u,\phi)$ where $p:C\otimes C\to C$ --{\em the product}-- and $u:\Bbbk\to
C$ --{\em the unit}-- are
coalgebra morphisms, and $\phi:C\otimes C\otimes C\to \Bbbk$ --{\em
  the associator}-- is a
convolution--invertible functional, satisfying the following axioms. 

\begin{equation}
  \label{e:unitprod}
  p(u\otimes \id)=\id=p(\id\otimes u)
\end{equation}
\begin{equation}\label{e:associator1} 
\sum (c_1 d_1)e_1 \phi(c_2 \otimes d_2 \otimes e_2)=\sum \phi(c_1
\otimes d_1 \otimes e_1)c_2(d_2e_2)
\end{equation} 
\begin{multline}\label{e:associator2}
\sum \phi(c_1d_1 \otimes e_1 \otimes f_1)\phi(c_2 \otimes d_2 \otimes
e_2f_2)=\\
= \sum \phi(c_1 \otimes d_1 \otimes e_1)\phi(c_2 \otimes d_2e_2
\otimes f_1)\phi(d_3 \otimes e_3 \otimes f_2)
\end{multline}
\begin{equation}\label{e:associator3}
\phi(c \otimes 1 \otimes d)=\varepsilon(c)\varepsilon(d)
\end{equation} 
The quadruple $(C,p,u,\phi)$ is called a {\em coquasi bialgebra}. 
\end{defi}

Along this paper coalgebras and coquasi bialgebras will be denoted
with the letters, $C$,$D$,\,etc.

The dual concept of a quasi bialgebra was originally defined in \cite{kn:dri}.

In the above equations we have written $1 \in C$ for the image under $u$
of the unit of $\Bbbk$, and $p(c,d)=cd$. Moreover,
when multiplying three elements of $C$ we used parenthesis
in order to establish the way we performed the operations.

The equation \eqref{e:associator1} can be interpreted in a precise
way as a naturality
condition on $\phi$. Equations \eqref{e:associator2} and
\eqref{e:associator3} can be written as the equalities $\phi(p\otimes
\id\otimes \id)\star\phi(\id\otimes\id\otimes
p)=(\phi\otimes\varepsilon)\star\phi(\id\otimes
p\otimes\id)\star(\varepsilon\otimes\phi)$  and 
$\phi(\id\otimes u\otimes 
\id)=\varepsilon\otimes\varepsilon $ valid in the convolution algebras
$(C\otimes C\otimes C)^\vee$ and $(C\otimes C)^\vee$ respectively.

Applying equation \eqref{e:associator2} to the case of $c=d=1$, 
we obtain that  $\phi(1 \otimes e \otimes f)= \phi(1 \otimes e_1 \otimes
f_1) \phi(1 \otimes e_2 \otimes f_2)$. Hence $\rho: C \otimes C
\rightarrow \Bbbk,\, \rho(e \otimes f) =
\phi(1 \otimes e \otimes f)$ is convolution invertible and $\rho
\star \rho = \rho$. Then $\rho=\varepsilon \otimes \varepsilon$ and
for later use we record below this and other similar consequence of
the axioms of a coquasi Hopf algebra.

 \begin{equation}\label{e:consequences}
 \phi(1 \otimes c \otimes d) = \varepsilon(c) \varepsilon(d)
 \hspace*{1cm} \phi(c \otimes d \otimes 1) = \varepsilon(c)
 \varepsilon(d)
 \end{equation}
 
\begin{obse} \label{obse:hcop}
If $(C,p,u,\phi)$ is a coquasi bialgebra, then
$C^{\mathrm{cop}}$ has a structure of a coquasi bialgebra with unit
$u$, multiplication $p \,\mathrm{sw}$ and associator
$\phi(\id \otimes
\mathrm{sw})(\mathrm{sw}\otimes\id)(\id\otimes\mathrm{sw}) $. 
We shall denote this coquasi bialgebra by $C^\circ$. In the literature
$C^\circ$ is denoted by $C^{\mathrm{copop}}$.
\end{obse}

Next we define the concept of {\em monoidal morphism} 
between coquasi bialgebras. 
Monoidal morphisms are to 
coquasi bialgebras what bialgebra morphisms
are to bialgebras.    

The monoidal morphisms are the adequate 
kind of morphisms for our category as they 
preserve multiplication and unit up to coherent isomorphisms. Although
we will only need this concept of morphism, for the sake of clarification 
we also give the definition of
{\em lax monoidal}\/ morphism as in this general case the role of the invertible scalar
$\rho$ appearing in the definition below is more transparent.

\begin{defi}\label{d:monstr}
Let $C$ and $D$ be coquasi bialgebras and $f:C\to D$ be a morphism of
coalgebras. A {\em lax monoidal structure}\/ on $f$ is a 
functional $\chi:C\otimes C\to\Bbbk$ and a
scalar $\rho \in \Bbbk$ satisfying
\begin{equation}
  \label{eq:chi1}
  (\chi\otimes p(f\otimes f) )\Delta_{C\otimes C}=(fp\otimes
  \chi)\Delta_{C\otimes C}  \qquad  u =  fu
\end{equation}
\begin{equation}  \label{eq:chimonoidal}
(\phi_D(f\otimes f\otimes
f))\star(\chi\otimes\varepsilon)\star(\chi(p\otimes \operatorname {id}))=
(\varepsilon\otimes\chi)\star(\chi(\id \otimes p))\star\phi_C
\end{equation}
\begin{equation}
  \label{eq:rhomonoidal}
 \rho \chi(u\otimes\id) = \varepsilon=\rho \chi(\id\otimes u).
\end{equation}
The lax monoidal structure $(\chi,\rho)$ is called a {\em monoidal
  structure}\/ when $\chi$ is invertible, -- notice that $\rho$ is
always invertible--. A morphism
of coalgebras between coquasi bialgebras equipped with a (lax)
monoidal structure is called a {\em (lax)  monoidal morphism}. 

A monoidal morphism between two coquasi bialgebras is a morphism of
coalgebras equipped with a monoidal structure. 
\end{defi}

The above terminology on monoidal structures comes
from category theory. In fact, the concept of 
monoidal morphism as defined above is a special
instance of a the concept of monoidal 1-cell between pseudomonoids.

In particular, 
in the case that the monoidal structure is $\chi=\varepsilon \otimes
\varepsilon$ and $\rho=1$ equations $\eqref{eq:chi1}$,
$\eqref{eq:chimonoidal}$ and $\eqref{eq:rhomonoidal}$ simply say that $f$
preserves the product, the unit and that $\phi_D(f \otimes f \otimes
f)=\phi_C$. 
Hence, it is clear in particular
that a map of coquasi bialgebras that preserves the product, the
coproduct and the associator, is a monoidal morphism.

The definition above can be generalized to the concept of 
a monoidal $(C,D)$--bicomodule, where $C$ and $D$
are coquasi bialgebras. In that case, 
$f:C\to D$ is a monoidal morphism if and only if the
bicomodule $f_+$ (see Definition
\ref{defi:f+}) is a monoidal bicomodule (see Theorem
\ref{theo:fmonoidal}). 
In this paper we will not cover these general aspects of the theory.
\begin{obse}\label{o:compmonmorph}
If $f:C\to D$ and $g:D\to E$ are monoidal morphisms with
monoidal structures $(\chi^f,\rho^f)$ and $(\chi^g,\rho^g)$
respectively, then $gf$ has canonical monoidal structure, namely, 
$(\chi^g(f\otimes f)\star\chi^f,\rho^f\rho^g)$. Also, the identity
morphism $\id:C\to C$ is equipped with a monoidal structure given by
$(\varepsilon\otimes\varepsilon,1)$. 
\end{obse}

In Proposition \ref{obs:dualsrlcirc2} we show that the antipode $S$
--see the definition below-- of a finite
dimensional coquasi Hopf algebra $H$ is a monoidal morphism from
$H^{\mathbf{copop}}$ to $H$. 
\begin{defi}
An {\em antipode}\/ for the coquasi bialgebra $H$ is a triple
$(S,\alpha,\beta)$ where $S: H^{\mathrm{cop}} \rightarrow H$ is a
  coalgebra morphism and the functionals $\alpha,\beta: H \rightarrow
  \Bbbk$ satisfy the following equations.

\begin{equation}\label{e:antipode1}
\sum S(h_1) \alpha(h_2) h_3=\alpha(h)1 \hspace*{1cm} \sum h_1\beta(h_2)S(h_3)=\beta(h)1
\end{equation}
\begin{equation}\label{e:antipode2}
\sum \phi^{-1}(h_1 \otimes Sh_3 \otimes h_5)\beta(h_2)\alpha(h_4)=\varepsilon(h)
\end{equation}
\begin{equation}\label{e:antipode3}
 \sum \phi(Sh_1 \otimes h_3 \otimes
 Sh_5)\alpha(h_2)\beta(h_4)=\varepsilon(h)
\end{equation}

A {\em coquasi Hopf algebra}\/ is a coquasi bialgebra equipped  with
an antipode.
\end{defi}
Along this paper coquasi Hopf algebras will be denoted as $H$.

\begin{obse} \label{obse:hcopant} If $(S,\alpha,\beta)$ is an antipode for the
  coquasi bialgebra $H$ then $(S,\beta,\alpha)$ is an antipode for
  the coquasi bialgebra $H^{\circ}$ considered in Observation 
\ref{obse:hcop}.
  
\end{obse}
\begin{obse}\label{obse:grouplike} For future use we record the
  following fact. If $a \in H$ is a group like element then
  $S(a)=a^{-1}$. Indeed, from the equality 
$\phi(Sa\otimes a\otimes Sa)\alpha(a)\beta(a)=1$, we deduce that 
$\alpha(a) \neq 0$. 
Then, from the equality
$S(a)\alpha(a)a=\alpha(a)1$ we deduce that $S(a)a=1$. The equality
$aS(a)=1$ can be proved in a similar manner.

Moreover, if $b \in H$ satisfies that $ab=ba=1$, then
$Sa=(ba)Sa$. Reassociating the product, 
if we call $\gamma:H \rightarrow \Bbbk$ the functional
$\gamma(x)= \phi(x \otimes a \otimes Sa)$, we have that:
$Sa=(ba)Sa=\gamma^{-1} \rightharpoonup b \leftharpoonup \gamma$. Then
$b=\gamma \rightharpoonup Sa \leftharpoonup \gamma^{-1}$ and  
since $Sa$ is a group like element, 
$b=\gamma  \rightharpoonup Sa \leftharpoonup \gamma^{-1}=Sa$. 

In the particular case of the group like element $1$, we have that
$S(1)=1$. 
\end{obse}

The definition of coquasi bialgebra (or rather its dual concept of
quasi bi\-al\-ge\-bra) was introduced by Drinfel'd in \cite{kn:dri}. The
crucial observation was that in order to guarantee the 
corresponding module category to be
monoidal, the associativity of the coproduct was only necessary up
to conjugation. The
concept of antipode (called by many authors a quasi antipode) 
and of quasi Hopf algebra as defined in
\cite{kn:dri}, is needed in order
guarantee the existence of duals in the corresponding categories
of finite dimensional objects. 

Along this paper we study the basic properties of the  
categories of modules and comodules for a coquasi Hopf algebra. These
categories have been considered by many authors
--see for example \cite{kn:kassel}
and more specifically \cite{kn:hn} and
\cite{kn:Schauenburg}.

Our main interest lays in the case that the
coquasi Hopf algebra is finite dimensional as a vector space. 
In this case, it is known (see
\cite{kn:bc} and \cite{kn:Schauenburg}) that the antipode $S$
is a bijective linear transformation. The composition inverse of $S$ will be
denoted as $\overline{S}$. 

We finish this Introduction by describing some of the notations we
use.  

We denote the usual duality functor in
the category of vector spaces as $V \mapsto V^{\vee}$ and the usual
evaluation and coevaluation maps as $\mathrm {e}$ and $\mathrm {c}$.   

Let $\mathcal C= (\mathcal C,  \otimes , \Bbbk, \Phi, l,r)$ 
be a monoidal category with monoidal structure $\otimes: \mathcal
C \times \mathcal C \rightarrow \mathcal C$, unit object
object $\Bbbk$, associativity constraint with components $\Phi_{M,N,L}: (M
\otimes N) \otimes L \rightarrow M \otimes (N \otimes L)$ 
and left and right unit
constraints $l$ and $r$. We denote as $\mathcal C^{\mathrm {rev}}$ the
tensor category $(\mathcal C, \otimes^{\mathrm{rev}} = \otimes
\mathrm{sw}, \Bbbk, \widehat{\Phi}, \widehat {r},\widehat{l})$ where
$\otimes^{\mathrm{rev}}(M \otimes N)= 
N \otimes M$, $\widehat{\Phi}_{M,N,L}= \Phi_{L,N,M}^{-1}$, $\widehat{r}=l$
and $\widehat{l}=r$. 

Let $\mathcal C$ and $\mathcal D$ be monoidal categories and
$T:\mathcal C\to\mathcal D$ a functor.
A monoidal structure on $T$ is a natural isomorphism $\otimes (T\times
T)\Rightarrow T\,\otimes:\mathcal C\times\mathcal C\to\mathcal D$ and an
isomorphism 
$\Bbbk \to T(\Bbbk)$ satisfying a certain natural list of
coherence axioms (see \cite{kn:JS} for details).
A monoidal functor is a functor equipped with a monoidal structure. 

We assume that the reader is familiar with the basic concepts
concerning rigidity for tensor categories as presented for example in
\cite{kn:JS} or \cite{kn:kassel}.
Recall that monoidal functors preserve duals. In other words,
if $T:\mathcal C\to\mathcal D$ is a monoidal functor and $M \in
\mathcal C$ is a left rigid object in $\mathcal C$, then $T(M)$ is
also left rigid and there is a
canonical natural isomorphism $\eta: T({}^*M)\to {}^*T(M)$. This isomorphism
is the unique arrow such that makes the diagram below commutative 
\[
\diagramcompileto{Tpresduals}
T({}^*M) \otimes T(M) \ar[d]_{\eta \otimes \operatorname{id}} \ar[r]^-{a}& 
T({}^*M \otimes M)\ar[r]^-{T(\operatorname{ev}_M)}& T(\Bbbk) \ar[d]^{b} \\ 
 {}^*T(M) \otimes
 T(M)\ar[rr]_-{\operatorname{ev}_{T(M)}}&&\Bbbk
\enddiagram
\]
where the maps $a$ and $b$ are the maps given by the monoidal
structure of $T$. 

Also, following the standard usage we write
${}^C\!\M$, $\M^D$  and ${}^C\!\mathcal M^D$ for the categories of left
$C$--comodules, right $D$--comodules and $(C,D)$--bicomodules, where $C$
and $D$ are coalgebras. For the coactions associated to the objects
in these categories we also use Sweedler's notation. In these
situations, when we add the subscript $f$ to the symbols, we mean to
say that we are
restricting our attention to the subcategories of finite--dimensional
objects, {\em e.g.}, ${}^C\!\M^D_f$ is
the full subcategory of ${}^C\!\M^D$ whose objects are finite
dimensional $\Bbbk$--spaces.  

\section{The categories of bicomodules and of Hopf modules}\label{s:catsofbicomod}
\subsection{The cotensor product}

We start by briefly reviewing some of the basic properties of
the cotensor product.
Given coalgebras $C,D$ and $E$, one can define the {\em cotensor
  product}\/ functor ${}^C\!\M^D\times{}^D\!\M^E\to{}^C\!\M^D$ as follows. 
If $M$ and $N$ are objects of ${}^C\!\M^D$ and ${}^D\!\M^E$ respectively,
its cotensor product over $D$, denoted by $M\square_DN$ is the
equalizer of the following diagram 
$${
\diagramcompileto{cotensor}
M\otimes N\ar@<3pt>[rrr]^-{((\varepsilon\otimes\id\otimes
  \id)\chi_M)\otimes\id} \ar@<-3pt>[rrr]_-{\id\otimes
  ((\id\otimes\id\otimes\varepsilon)\chi_N)}&&&
M\otimes D\otimes N
\enddiagram}
$$
endowed with the bicomodule structure induced by the left coaction of
$M$ and the right coaction of $N$. 

If $F$ is another coalgebra, there is a natural isomorphism between
the two obvious functors
${}^C\!\M^D\times{}^D\!\M^E\times{}^E\!\M^F\to{}^C\!\M^F$, with components
$L\square_D(M\square_EN)\cong(L\square_DM)\square_EN$ induced by the
universal property of the equalizers. Also, the functors
$C\square_C-$ and $-\square_DD:{}^C\!\M^D\to{}^C\!\M^D$ are canonically
isomorphic to the identity functor. All these data satisfy coherence
conditions\,; in categorical terminology we say that the categories of bicomodules
form a {\em bicategory}\/ \cite{kn:Benabou}.

From the above, it is clear that the cotensor product provides a
monoidal structure to ${}^C\!\M^C$, with unit object $(C,\Delta^2)$. 

\begin{obse}\label{o:cotensorfinitary}
  The cotensor product functor
  ${}^C\!\M^D\times{}^D\!\M^E\to{}^C\!\M^D$ preserves filtered
  colimits in each variable. This is because finite limits commute with
  finite colimits. 
\end{obse}
Next we consider corestriction functors.

 \begin{defi} \label{defi:f+}
If $f:C\to D$ is a morphism of coalgebras, we shall denote by $f_+ =
 C_f  \in
 {}^C\!\M^D$
 the object obtained from the regular bicomodule $C$ by
 correstriction with $f$ on the right, {\em i.e.}, 
the coaction in $f_+=C_f$ is given
 by $x\mapsto \sum x_1\otimes x_2\otimes f(x_3)$. Similarly, we
 shall denote by $f^+ ={}_f C \in {}^D\!\M^C$ the object obtained form $C$ by
 correstriction on the left. 
\end{defi}
 
 Taking cotensor products with bicomodules that are induced by morphisms of
 coalgebras has convenient properties. 
 \begin{obse}\label{obse:M_f}
Suppose that $f, C$ and $D$ are as above and that 
$A$ is an arbitrary coalgebra,
\begin{enumerate}
\item  For any bicomodule $M\in{}^A\!\M^C$, with coaction $\chi_M$, the
   cotensor product $M\square_C f_+ = M\square_C C_f 
\in {}^A\!\M^D$ is canonically
   isomorphic with the
   bicomodule --sometimes called also $M_f$-- with underlying space 
   $M$ and coaction
   $(\id_A\otimes\id_M\otimes f)\chi_M:M\to A\otimes M\otimes D$. 

 In a completely analogous way, if $N\in {}^C\!\M^A$, then the cotensor
 product $f^+\square _C N = {}_f C \square_C N \in {}^D\!\M^A$ is canonically isomorphic to the
 bicomodule --sometimes called ${}_f N$-- 
with underlying space $N$ and coaction
 $(f\otimes\id_N\otimes \id_C)\chi_N:N\to D\otimes N\otimes A$.

 Hence, $-\square_Cf_+ = -\square_C C_f $ and $f^+\square_C- = {}_f C
 \square_C-$ are the functors $M \mapsto M_f: {}^A\!\M^C \rightarrow
 {}^A\!\M^D$ and  $N \mapsto {}_fN: {}^C\!\M^A \rightarrow {}^D\!\M^A$
 given by correstriction with $f$. 

\item Given a morphism of coalgebras $f:C\to D$, it is clear that 
 the functor considered above $-\square_C
   f_+ = -\square_C C_f
   :{}^A\M^C\to{}^A\M^D$ is left adjoint to the so called coinduction
   functor $-\square_D f^+ =
   -\square_D {}_f C : {}^A\M^D
   \rightarrow {}^A\M^C$.

Similarly, the other correstriction functor   $f^+\square _C - = {}_f
C \square_C -: {}^C\!\M^A \rightarrow  {}^D\!\M^A$ is left adjoint to
the so called coinduction functor $f_+ \square _D - = C_f 
\square_C -: {}^D\!\M^A \rightarrow  {}^C\!\M^A$
\item Assume that we have two morphisms of coalgebras $f:C \rightarrow
  D$ and $g:D \rightarrow E$. In that situation we have canonical
  isomorphisms between $(gf)_+ \cong f_+ \square_D g_+$ and $(gf)^+
  \cong g^+ \square_D f^+$.
\end{enumerate}
 \end{obse}

 \begin{defi}
 Assume that
 the coalgebra $D$ has a group like element that we call $1 \in D$. We
 apply the above construction to the morphism of coalgebras $u: \Bbbk
 \rightarrow D$. In this case we abbreviate $(-)_0=-\square_{\Bbbk}
 u_+=(-)_u: {}^A\!\M \rightarrow {}^A\!\M^D$. Similarly we call 
${}_0(-)=u^+ \square_{\Bbbk}- ={}_u(-): \M^A
 \rightarrow {}^D\!\M^A$. 
 \end{defi}

\begin{obse}
\begin{enumerate}
\item The underlying space for $M_0$ is $M$ and the explicit formula for the 
 associated coaction is $\chi_0(m)=\sum m_{-1}\otimes m_0 \otimes 1$ if
 $(M,\chi) \in {}^A\!\M$. As we observed before this functor is left
 adjoint to $-\square_{\Bbbk} u^+: {}^A\!\M^D \rightarrow {}^A\!\M$. 
 It is clear that if $M \in {}^A\!\M^D $, then $N\square_{\Bbbk} u^+=
 N^{{\mathrm{co}}D}$. In other words, the functor $(-)_0$ that produces
 from a left $A$--comodule the $(A,D)$--bicomodule with trivial right
 structure is left adjoint to the fixed point functor. 
\item The underlying space to ${}_0M$ is the same
 than $M$ and the coaction is   
is ${}_0\chi (m)=\sum 1 \otimes m_{0} \otimes m_1$.
This functor is left
 adjoint to $u_+ \square_{\Bbbk} -: {}^D\!\M^A \rightarrow \M^A$, 
that is the functor that takes the left coinvariants,
 i.e., sends $M \mapsto {}^{\mathrm{co}D}M$. 
 In other words, the functor $(-)_0$ that produces
 from a right $A$--comodule the $(D,A)$--bicomodule with the trivial left
 structure is left adjoint to the left fixed point functor. 
\end{enumerate}
\end{obse}

\begin{theo} \label{l:nacho}
Let $g,h:C \rightarrow D$ be two morphisms of 
  coalgebras, and consider the following structures.  
\begin{enumerate}
\item \label{l:nachoit1}
Functionals $\gamma:C\to\Bbbk$ satisfying
$(\gamma\otimes g)\Delta=(h\otimes\gamma)\Delta$.
\item \label{l:nachoit2}
 Morphisms of bicomodules $\theta:g_+ \to h_+ $.
\item \label{l:nachoit3} 
Natural transformations
  $\Theta:(-\square_Cg_+)\Rightarrow(-\square_Ch_+):\M^C\to\M^D$.
\end{enumerate}
Each structure of type \eqref{l:nachoit1} induces a structure
of type \eqref{l:nachoit2} by $\theta =  (\mathrm{id}_C \otimes \gamma)\Delta$
and each structure of type \eqref{l:nachoit2} induces
a structure of type \eqref{l:nachoit3} by $\Theta
  = -\square_C \theta$. 
Moreover, if $C$ is finite dimensional 
  these correspondences are bijections, with inverses given by
  $\gamma=\varepsilon\theta$ and $\theta=\Theta_C$. 
The identity and the composition of the natural transformations in
\eqref{l:nachoit3} correspond
to the identity and the composition of the morphisms in
\eqref{l:nachoit2} and to the counit
$\varepsilon$ and  the convolution product of
the functionals in \eqref{l:nachoit1}. 
In particular,  the natural transformation $\Theta$ is invertible iff the
 associated 
 morphism $\theta$ is
invertible iff the corresponding functional $\gamma$ is
convolution--invertible. 
\end{theo}
\begin{proof}
That each structure \eqref{l:nachoit1} induces a structure
\eqref{l:nachoit2} and each structure in \eqref{l:nachoit2} induces
a structure \eqref{l:nachoit3} is easily verified. 

Next we prove that the above described maps are indeed a 
bijection between \eqref{l:nachoit1} and
\eqref{l:nachoit2}. 
The map $\theta $ satisfies 
\begin{equation}\label{eq:thetamor}
(\id\otimes\theta\otimes
g)\Delta^2=(\id\otimes\id\otimes h)\Delta^2\theta.
\end{equation} Composing the
above equality with $\id \otimes \id \otimes \varepsilon$ we deduce
that $\Delta \theta= (\id \otimes \theta) \Delta$. Now, if we compose
Equation \eqref{eq:thetamor} with 
$\varepsilon \otimes \varepsilon \otimes \id$
we obtain
$(\gamma\otimes g)\Delta=h\theta$ , with $\gamma = \varepsilon
\theta$. A direct calculation shows that 
$(h\otimes\gamma)\Delta= (h\otimes\varepsilon)(\id
\otimes \theta)\Delta= (h\otimes\varepsilon)\Delta \theta= h\theta$.
Hence, $\gamma=\varepsilon \theta$ satisfies condition (1). Clearly,
the correspondences given above between elements $\gamma$ and $\theta$
are inverses of each other.

If we assume that $C$ is finite dimensional, 
the bijection between the structures in \eqref{l:nachoit2} and
\eqref{l:nachoit3} is a consequence of Observation \ref{o:apendix}. 
\end{proof}

For a functional $\gamma$ as in Theorem
\ref{l:nacho}.3, we will denote as $\gamma_+:g_+\to
h_+$ the associated morphism
of comodules and as $\Gamma$ the corresponding natural transformation.  

\begin{obse}\label{obse:concrete} In the case that $\gamma$ is
  convolution invertible, the condition (3) that relates  $g$ and
  $h$ in Theorem \ref{l:nacho}, can be written as $g = \gamma^{-1}\star h
  \star \gamma$\, or as any of the equalities below valid for all $c\in C$:  
\begin{equation}\label{e:concrete} g(c \leftharpoonup \gamma)= h(\gamma
  \rightharpoonup c) \quad \/,\/ \hspace*{.5cm}g(c)= h(\gamma
  \rightharpoonup c \leftharpoonup \gamma^{-1}) 
\end{equation}
\end{obse}

\subsection{The tensor product over $\Bbbk$}
When we consider  coalgebras that have the additional structure of a 
coquasi bialgebra, 
the corresponding 
categories of comodules and of bicomodules have 
--besides the monoidal structure given by the cotensor product--
another monoidal structure.  
This monoidal structure is based  upon the tensor 
product over the base field $\Bbbk$ with associativity constraint
defined in terms of the corresponding functional $\phi$. 
For example if $C$ and $D$ are coquasi bialgebras with associators
$\phi_C$ and $\phi_D$ respectively, if $L,M,N \in {}^C\!\M^D$ the
associativity constraint is the map 
\begin{equation*}
\Phi_{L,M,N}: (L \otimes M)
\otimes N \rightarrow L \otimes (M \otimes N) 
\end{equation*} 
given by the formula 
\begin{equation}\label{e:constraint}
\Phi((l \otimes m)
\otimes n)= \sum \phi_C(l_{-1} \otimes m_{-1} \otimes n_{-1}) l_0
\otimes (m_0 \otimes n_0) \phi_D^{-1}(l_{1} \otimes m_{1} \otimes n_1)
\end{equation}  
Here we view $M \otimes N$ as an
object in ${}^C\!\M^D$ with the usual structure: $\chi_{M \otimes
  N}(m \otimes n)= \sum m_{-1}n_{-1} \otimes m_0 \otimes n_0 \otimes m_1n_1
\in C \otimes M \otimes N \otimes D$. Notice that the above formula
for the associativity constraint can be written using the standard 
actions associated to coactions as follows: 
\begin{equation*}\Phi_{L,M,N}((l \otimes m) \otimes n)= \phi_D^{-1}
\rightharpoonup l \otimes m \otimes n \leftharpoonup \phi_C \in L
\otimes (M \otimes N).
\end{equation*} 

The unit constraints $M\otimes \Bbbk \cong M$ and $\Bbbk\otimes M\cong
M$ are the same than in the category of $\Bbbk$--vector spaces. 

In case that the categories are ${}^C\!\M$ or $\M^C$, the
constraints are defined similarly but using only the action by $\phi^{-1}$
on the left for $\M^C$ and of $\phi$ on the right for ${}^C\!\M$. 
The monoidal categories
$\M^C$ and ${}^C\!\M$ can also be defined as 
${}^\Bbbk\!\M^C$ and ${}^C\!\M^\Bbbk$ respectively.

\begin{obse}\label{obse:associativity} If $C$ is a coquasi bialgebra
  with unit $u$ and multiplication $p$, the triple $(C,p,u) \in 
  {}^C\mathcal M\,^C$ is an associative algebra. This can be proved directly using
  equation (\ref{e:associator1}), which can be rewritten as  $p(\id
  \otimes p)\Phi_{C,C,C}= p(p \otimes \id): (C \otimes C) \otimes
    C \rightarrow C$. 
\end{obse}

\begin{defi}\label{defi:hopfmod}
The category of right $C$--modules within $\bicomod$ will be denoted
as $\hopfmod$ and called the category of Hopf modules. Similarly we
define the category ${}^C_C \mathcal M^C$.
\end{defi}

Notice that the unit object $\Bbbk$ of the monoidal structure $\otimes$
is canonically a Hopf module with action given by the counit $\varepsilon$.

The category of Hopf modules in this context was first considered 
in \cite{kn:hn} and \cite{kn:Schauenburg}.

\begin{obse}\label{obse:squareandtensor}
\noindent 
a) The $\square_C$ monoidal structure of $\bicomod$ lifts to a monoidal
structure on $\hopfmod$ in such a way that the forgetful functor
$\hopfmod\to\bicomod$ is monoidal. 

Indeed, if $M,N,L,R$ are in $\bicomod$,
one easily can define --using the universal property of equalizers-- 
a natural morphism of bicomodules $(M\,\square_C N)\otimes
(L\,\square_C R)\to (M\otimes L)\,\square_C(N\otimes R)$ relating both
monoidal structures on $\bicomod$. If $L=R=C$, we obtain a map
$
(M\,\square_CN)\otimes C\to (M\otimes C)\,\square_C(N\otimes C)
$ that
composed with the right $C$--actions on $M$ and $N$ ---$a_M: M \otimes C
\rightarrow M$ and $a_N: N \otimes C \rightarrow N$---endows $M\,\square_C
N$ with the structure of a Hopf module
\begin{equation} (M\,\square_CN)\otimes C\to (M\otimes
  C)\,\square_C(N\otimes C) \xrightarrow{a_M \otimes a_N} M\,\square_CN. 
\end{equation} Clearly the unit object $C$ of
$\square_C$ in $\bicomod$ is also a unit object in $\hopfmod$.

\noindent
b) If $f: \Bbbk \rightarrow C$ and $g: C \rightarrow C$ are morphisms
of coalgebras --in particular this means that $f(1) \in C$ is a group
like element-- then $f_+ \otimes (M \/\square_C\, g_+) \cong M \/\square_C\,
p(f \otimes g)_+$ and $(M \square_C\, g_+) \otimes f_+ \cong M \square_C\,
p(g \otimes f)_+$.
\end{obse}

\subsection{Monoidal functors induced by monoidal morphisms}
\label{ss:monmods}
\begin{theo}\label{theo:fmonoidal}
  Let $f:C\to D$ be a coalgebra morphism, and consider the following
  structures.
\begin{enumerate}
  \item Monoidal structures on $f$,
   \item Monoidal structures on the functor $(-\square_C
      f_+):\mathcal M^C\to\mathcal M^D$,
   \item Monoidal structures on the functor
     $(f^+\square_C-):{}^C\mathcal M\to{}^D\mathcal M$.
  \end{enumerate}
Each structure (1) induces structures (2) and (3). Moreover, if
$C$ is finite-dimensional there is a bijection between the three types
of structures.
\end{theo}
\begin{proof}
First, we consider the relationship of structures of type (1) with
structures of type (2). 
If $\chi:C\otimes C\to\Bbbk$, $\rho\in\Bbbk$ is a monoidal structure
on $f:C\to D$, then the transformation with components
$\Theta_{M,N}:M\otimes N\to M\otimes N$ given by $\Theta_{M,N}(m
\otimes n) =\chi \rightharpoonup (m \otimes n) = \sum \chi(m_1\otimes n_1)
m_0\otimes n_0$ together with the isomorphism
$\Bbbk\to\Bbbk$ given by multiplication by $\rho$ is a monoidal
structure as in (2). Indeed, for example the condition that 
the map $\Theta_{M,N}:M_f \otimes N_f
\rightarrow (M \otimes N)_f$ is a morphism of
$H$--comodules --recall the notations of Observation \ref{obse:M_f}--
is equivalent with condition \eqref{eq:chi1} in Definition
\ref{d:monstr}.   
A structure as in (3) is obtained in a similar
way, using $\chi^{-1}$ and $\rho^{-1}$. 

If we now assume that $C$ is finite-dimensional we can proceed
backwards in order to go from (2) to (1). Let 
$\Theta$ be a natural
transformation as depicted below.
\begin{equation}\label{eq:Theta}
\xymatrixcolsep{1.3cm}
\diagramcompileto{monstr1}
{\M}^C\times{\M}^C\ar[rr]^-{(-\square_Cf_+)\times(-\square_Cf_+)}\ar[d]_\otimes
\drrtwocell<\omit>{\Theta}&&
{\M}^D\times{\M}^D\ar[d]^\otimes\\
{\M}^C\ar[rr]_-{(-\square_Cf_+)}&&
{\M}^D
\enddiagram
\end{equation}
Since all the functors in this diagram preserve filtered colimits,
$\Theta$ is determined by its restriction to the categories of
finite-dimensional comodules. So we can substitute the categories of
comodules in diagram \eqref{eq:Theta} by the corresponding categories of
finite-dimensional comodules. In the appendix --Section \ref{s:back}--
we prove that for
a finite dimensional coalgebra $C$, 
composition with the tensor product functor
$\otimes_\Bbbk:\M^C_f\times\M^C_f\to\M^{C\otimes C}_f$ induces an
equivalence $\mathrm{Lex}[\M^{C\otimes C}_f,\M^D_f]\simeq
\mathrm{Lex}[\M^C_f,\M^C_f;\M^D_f]$. The category 
$\mathrm{Lex}[\M^C_f,\M^C_f;\M^D_f]$ appearing on the
right hand side of the equivalence is the category of functors from
$\M^C_f\otimes\M^C_f$ to $\M^D_f$ which are left exact in each
variable--see also the Appendix for the definition of
$\M^C_f\otimes\M^C_f$--. 
Using this fact, we deduce that natural transformations as
in diagram \eqref{eq:Theta} are in bijection with natural
transformations as in the diagram below: 
\begin{equation*}
\xymatrixcolsep{1.2cm}
\diagramcompileto{Theta2}
\M^{C\otimes C}_f\ar[rr]^-{-\square_{C\otimes C}(f\otimes f)_+}
  \ar[d]_{-\square_{C\otimes C}p_+}\drrtwocell<\omit>{\Theta'}&&
\M^{D\otimes D}_f \ar[d]^{-\square_{D\otimes D}p_+} \\
\M^C_f\ar[rr]_-{-\square_C f_+}&&
\M^D_f.
\enddiagram
\end{equation*}
Also these type of natural transformations are in bijective
correspondence with bicomodule morphisms $(p(f\otimes f))_+\to
(fp)_+$. These bicomodule morphisms correspond bijectively  
to functionals $\chi:C\otimes C\to\Bbbk$ satisfying $\sum
\chi(c_1\otimes c'_1) f(c_2)f(c'_2)=\sum f(c_1c'_1)\chi(c_2\otimes
c'_2)$. The invertibility of $\Theta$ is equivalent to the
invertibility of $\chi$. 

Similarly, an isomorphism $\Sigma:\Bbbk\to \Bbbk\/\,\,\square_C f_+$ 
is just an invertible
scalar $\rho$ such that $\rho f(1)=\rho 1$. The axioms of a monoidal
structure for $\Theta,\Sigma$ translate to the axioms of a
monoidal structure on $f$ for $\chi,\rho$. 

The relationship between the structures (1) and (3) is as follows. A
monoidal structure $\chi:C\otimes C\to\Bbbk,\rho\in\Bbbk$ on $f$
induces a monoidal structure on $(f^+\square_C-)$ given by 
the $D$-comodule morphism
$m\otimes n\mapsto \sum \chi^{-1}(m_{-1}\otimes n_{-1})m_0\otimes n_0: 
(f^+\square_C M)\otimes (f^+\square_C N)\to f^+\square_C(M\otimes N)$
and by $\lambda\mapsto \rho^{-1}\lambda: \Bbbk\to f^+\square_C\Bbbk$. 
The proof of
the converse, i.e. that 
when $C$ is finite-dimensional every monoidal structure on
$(f^+\square_C)$ is of this form for a unique $(\chi,\rho)$ is similar
to the one presented above for the case of right comodules. 
\end{proof}

\begin{coro}\label{coro:()_0monoidal}  In the situation above the
  functors $u^+ \square_{\Bbbk} - = {}_0\,(-) : \M^D \rightarrow
{}^C\!\M^D$, $(-)_0=-\square_{\Bbbk} u_+: {}^C\!\M
\rightarrow {}^C\!\M^D$, are monoidal. 
\end{coro} 
\begin{proof} It follows immediately from the fact that $u:\Bbbk
  \rightarrow C$ preserves the product, the unit and the associator. 
Explicitly, $u$ has a monoidal
  structure provided by $\chi=\mathrm{id}:\Bbbk\to\Bbbk$ and
  $\rho=1\in\Bbbk$. In this case, equalities \eqref{eq:chi1} and
  \eqref{eq:rhomonoidal} are trivial while condition
  \eqref{eq:chimonoidal} reads as
  $\phi(1\otimes1\otimes 1) =1$, which is true by \eqref{e:consequences}.
\end{proof}

\subsection{Some useful monoidal functors on comodule categories}
In this subsection we describe the functors we use along this work. 

\begin{defi}\label{defi:()circrl}
If $C$ and $D$ are coalgebras we define the functors
\begin{equation*}(-)^\circ:{}^C\M\,^D\to{}^{D\,^{\mathrm{cop}}}\!\M\,^{C\,^{\mathrm{cop}}}
\qquad (-)^r:(\M\,^C_f)^{\mathrm{op}}\to {}^C\!\M_f
\end{equation*}
\begin{equation*}
(-)^\ell:({}^C\!\M_f)^{\mathrm{op}}\to\M\,^C_f
\end{equation*}
The functor $(-)^\circ$ is the identity on arrows, and if
$M\in{}^C\!\M ^D$ with coaction $\chi(m)=\sum m_{-1}\otimes m_0\otimes m_1$,
then $M^\circ$ has $M$ as underlying space and coaction
$\chi^\circ(m)=\sum m_1\otimes m_0\otimes m_{-1}$. In the case when
$C,D$ are coquasi bialgebras, $(-)^\circ$ has a canonical structure of
a monoidal functor
$({}^C\!\M\,^D)^{\mathrm{rev}}\to{}^{D^{\circ}}\!\M\,^{{C^{\circ}}}$
given by the usual symmetry of vector spaces
$\mathrm{sw}:M^\circ\otimes N^\circ \cong(N\otimes M)^\circ$ and the
identity $\Bbbk\to\Bbbk^\circ$. 

The functor $(-)^r$ is defined as follows. If $M\in \M^C_f$, the
underlying space of $M^r$ is $M^\vee$, the linear dual of $M$. If
$\mathrm{c}$ and $\mathrm {e}$ denote the standard coevaluation and
evaluation, the
coaction for $M^r$ is: 
\begin{equation}\label{eq:(-)^r}
M^\vee\xrightarrow{\id\otimes \mathrm{c}}M^\vee\otimes M\otimes M^\vee
\xrightarrow{\id\otimes \chi\otimes\id}M^\vee\otimes M\otimes
C\otimes M^\vee\xrightarrow{\mathrm{e}\otimes \id\otimes \id}C\otimes M^\vee.
\end{equation}
On arrows, $(-)^r$ is given by the usual (linear) duality functor. 
We call $M^r$ the {\em right adjoint} of $M$. 
When $C$ is a coquasi bialgebra, $(-)^r$ has the
following canonical structure of a monoidal functor
$(\M^C_f)^{\mathrm{op}}\to{}^C\!\M_f$. 
The unit constraint is the
canonical isomorphism $\Bbbk\cong \Bbbk^\vee$; if $M,N\in\M^C_f$, then the
transformation $M^r\otimes N^r\to (M\otimes N)^r$ is given by the canonical
arrows $M^\vee\otimes N^\vee\to (M\otimes N)^\vee$, which are
isomorphisms by dimension considerations. We should remark that here we
are not thinking $M^\vee$ as a categorical dual of the vector space
$M$ but rather as the internal hom $\mathbf{Vect}(M,\Bbbk)$. This is
the reason why $(-)^r$ does not reverse the order of the tensor
products. 

 The definition
of $(-)^\ell$ is analogous, if $N \in {}^C\!\M_f$, then:
\begin{equation}\label{eq:(-)^l}
N^\vee\xrightarrow{\mathrm{c}\otimes \id}N^\vee\otimes N\otimes N^\vee
\xrightarrow{\id\otimes \chi\otimes\id}N^\vee\otimes C\otimes
N\otimes N^\vee\xrightarrow{\id\otimes\id\otimes\mathrm{e}}C\otimes N^\vee.
\end{equation}
 If $N\in{}^C\!\M_f$, we call $N^\ell$ the
{\em left adjoint}\/ of $N$. When $C$ is a coquasi bialgebra we have a
monoidal functor $(-)^\ell:({}^C\!\M_f)^{\mathrm{op}}\to\M^C_f$. 
\end{defi}

For future reference we record the following results that can be
proved directly. 

\begin{lema}
Observe that $(-)^r$ and $(-)^\ell$ are inverse monoidal equivalences and that
$(-)^{r\ell} = (-)^{\ell r}$.  The monoidal 
isomorphisms $M^{r\ell}\cong M\cong M^{\ell r}$ are just the canonical
linear isomorphisms $M\cong M^{\vee\vee}$.   
\end{lema}
\begin{lema}\label{l:relations}
  For any morphism of coalgebras $f:C\to D$, the diagrams in Figure
  \ref{fig:relations} 
  commute. If moreover $f$ is a monoidal morphism, the diagrams
  commute as diagrams of monoidal functors. 
\begin{figure}\label{fig:relations}
$${
\xymatrixcolsep{1.6cm}
\diagramcompileto{f_+andf^+}
(\M^C_f)^{\mathrm{op}}\ar[r]^-{(-\square_C
  f_+)^{\mathrm{op}}}\ar[d]_{(-)^r}&
(\M^D_f)^{\mathrm{op}}\ar[d]^{(-)^r}\\
{}^C\!\M_f\ar[r]_-{f^+\square_C-}&
{}^D\!\M_f
\enddiagram
}
\qquad
{\xymatrixcolsep{1.6cm}
\diagramcompileto{f_+andf^+2}
({}^C\!\M_f)^{\mathrm{op}}
\ar[r]^-{(f^+\square_C-)^{\mathrm{op}}}\ar[d]_{(-)^\ell}&
({}^D\!\M_f)^{\mathrm{op}}\ar[d]^{(-)^\ell}\\
\M^C_f\ar[r]_-{-\square_Cf_+}&
\M^D_f
\enddiagram
}
$$
$$
{\xymatrixcolsep{1.8cm}
\diagramcompileto{f_+andcirc}
(\M^C)^{\mathrm{rev}}\ar[r]^-{(-\square_Df_+)^{\mathrm{rev}}}
\ar[d]_{(-)^\circ}&
(\M^D)^{\mathrm{rev}}\ar[d]^{(-)^\circ}\\
{}^{C^{\circ}}\!\M\ar[r]_-{{f^{\mathrm{cop}}}^+\square_{C^{\mathrm{cop}}}-}&
{}^{D^{\circ}}\!\M
 \enddiagram}
$$
\caption{Diagrams in Lemma \ref{l:relations}}
\end{figure}
\end{lema}
\begin{proof}
  Recall that if $f$ has a monoidal structure $\chi:C\otimes
  C\to\Bbbk, \rho\in\Bbbk$, the monoidal structures on
  $(-\square_Cf_+)$ and $(f^+\square_C-)$ are induced by $\chi,\rho$
  and $\chi^{-1},\rho^{-1}$, respectively. Also, it is easy to show
  that the monoidal
  structures on $(-\square_Cf_+)^{\mathrm{op}}$ and
  $(f^+\square_C-)^{\mathrm{op}}$ are induced by $\chi^{-1},\rho^{-1}$
  and $\chi,\rho$ respectively. The
  verification of the Lemma is direct.
\end{proof}
\begin{lema}
  For any coquasi bialgebra $C$ the following diagram of monoidal
  functors commutes. 
$$
\xymatrixcolsep{1.2cm}
\diagramcompileto{circandrandell}
({}^C\!\M_f)^{\mathrm{oprev}}\ar[r]^-{((-)^\circ)^{\mathrm{op}}}
\ar[d]_{((-)^\ell)^{\mathrm{rev}}}&
(\M_f^{C^{\circ}})^{\mathrm{op}}\ar[d]^{(-)^r}\\
(\M^C_f)^{\mathrm{rev}}\ar[r]_-{(-)^\circ}&
{}^{C^{\circ}}\!\M_f
\enddiagram
$$
\end{lema}
\begin{lema}\label{l:roro}
  The following two monoidal functors are monoidally
  isomorphic to the identity functor via the canonical maps
$M \mapsto M^{\vee\vee}$ 
$$
\M^C_f\xrightarrow{(-)^r}{}^{C}\!\M_f^{\mathrm{op}}
\xrightarrow{(-)^\circ} (\M^{C^{\circ}}_f)^{\mathrm{oprev}}
\xrightarrow{(-)^r}{}^{C^{\circ}}\!\M_f^{\mathrm{rev}}
\xrightarrow{(-)^\circ}\M^C_f
$$
$$
{}^C\!\M_f\xrightarrow{(-)^\ell}(\M^C_f)^{\mathrm{op}}
\xrightarrow{(-)^\circ} ({}^{C^{\circ}}\!\M_f)^{\mathrm{oprev}}
\xrightarrow{(-)^\ell}(\M^{C^{\circ}}_f)^{\mathrm{rev}}
\xrightarrow{(-)^\circ}\M^C_f
$$
\end{lema}

\section{Duality}\label{s:duality}
In the case that the map $S$ is invertible --for example if the
coquasi Hopf algebra $H$ is
finite dimensional-- the monoidal categories $\lcomodH_f$, $\rcomodH_f$ 
and $\bicomodH_f$ are
rigid. In the case of $\bicomodH$, {\em e.g.}, we need to construct for
every object $M$ a left  and a right dual --denoted as
${}^*M$ and $M^*$ respectively-- together with the corresponding
evaluation and coevaluation maps. 

If $(M,\chi) \in \bicomodH_f$ and  $M^\vee$ is the dual of the
underlying vector space, $M^*=(M^\vee,\chi^*)$ where $\chi^*$ is the
composition 
\begin{equation}\label{e:M*structure}
\begin{split}&M^\vee \xrightarrow{\mathrm {c} \otimes\id}M^\vee\otimes M\otimes
M^\vee \xrightarrow{\id\otimes\chi\otimes\id} M^\vee\otimes H\otimes
M\otimes H \otimes M^\vee \xrightarrow{\mathrm{sw}\otimes\id\otimes
  \mathrm{sw}}\\&H\otimes
M^\vee\otimes M\otimes M^\vee \otimes H
\xrightarrow{\id\otimes\id\otimes \mathrm{e} \otimes \id}H\otimes
M^\vee\otimes H\xrightarrow{S\otimes\id\otimes \overline S}H\otimes
M^\vee\otimes H.
\end{split}
\end{equation}
The evaluation and coevaluation morphisms are given by
\begin{equation}
  \label{eq:ev}
  \mathrm{ev}^\ell:{}^*M\otimes M\xrightarrow{\id\otimes\chi}{}^*M\otimes
  H\otimes M\otimes H \xrightarrow{\id\otimes \beta\overline S\otimes \id\otimes
    \alpha} {}^*M\otimes M\xrightarrow{\operatorname{e}}\Bbbk
\end{equation}
and
\begin{equation}
  \label{eq:coev}
\mathrm{coev}^\ell:  \Bbbk\xrightarrow{\mathrm{c}}M\otimes{}^*M\xrightarrow{\chi\otimes
    \id}H\otimes M\otimes H\otimes {}^*M\xrightarrow{\alpha
\overline    S\otimes\id\otimes \beta\otimes\id}M\otimes{}^*M
\end{equation}
It is not hard to check using
\eqref{e:antipode1} that $\mathrm{ev}^\ell$ and
$\mathrm{coev}^\ell$ are morphisms in $\bicomodH$.
Moreover,
the maps $\mathrm{ev}^{\ell}
:{}^{*}M \otimes M \rightarrow \Bbbk \in \bicomodH$
and
$\mathrm{coev}^{\ell}: \Bbbk \rightarrow M \otimes {}^{*}M \in
\bicomodH$ satisfy 
\begin{equation}\label{eq:dualcond1}\id_{M}=
M\xrightarrow{\mathrm{coev}^{\ell}\otimes \id}(M \otimes {}^{*}M)
\otimes M\xrightarrow{\Phi_{M,{}^{*}M,M}}M\otimes({}^{*}M \otimes
M)\xrightarrow{\id\otimes\mathrm{ev}^{\ell}}M
\end{equation} 
and
\begin{equation}\label{eq:dualcond2}
\id_{{}^{*}M}=
{}^{*}M\xrightarrow{\id \otimes \mathrm{coev}^{\ell}}{}^{*}M\otimes(M
\otimes {}^{*}M)\xrightarrow{\Phi^{-1}_{{}^{*}M,M,{}^{*}M}}({}^{*}M \otimes
M)\otimes{}^{*}M\xrightarrow{\mathrm{ev}^{\ell}\otimes \id}{}^{*}M.
\end{equation}
These equations are direct consequences of \eqref{e:antipode2} and
\eqref{e:antipode3}. 

Analogously, ${}^*M=(M^\vee,{}^*\chi)$, where ${}^*\chi$ is the
composition
\begin{equation}\label{e:*Mstructure}
\begin{split}&M^\vee \xrightarrow{\id\otimes \mathrm{c}}M^\vee\otimes M\otimes
M^\vee \xrightarrow{\id\otimes\chi\otimes\id} M^\vee\otimes H\otimes
M\otimes H \otimes M^\vee
\xrightarrow{\mathrm{sw}\otimes\id\otimes\mathrm{sw}}\\ & H\otimes
M^\vee\otimes M\otimes M^\vee \otimes H
\xrightarrow{\id\otimes \mathrm{e}\otimes\id\otimes \id}H\otimes
M^\vee\otimes H\xrightarrow{\overline S\otimes\id\otimes  S}H\otimes
M^\vee\otimes H.
\end{split}
\end{equation}
The corresponding right evaluation and
coevaluation morphisms are:
\begin{equation}
  \label{eq:evr}
  \mathrm{ev}^r:M\otimes M^*\xrightarrow{\chi\otimes\id}
  H\otimes M\otimes H\otimes M^* \xrightarrow{\beta \otimes \id\otimes
    \alpha \overline S\otimes\id} M\otimes M^*\xrightarrow{\mathrm{e}}\Bbbk
\end{equation}
and
\begin{equation}
  \label{eq:coevr}
\mathrm{coev}^r:  \Bbbk\xrightarrow{\mathrm{c}}M^*\otimes
M\xrightarrow{\id\otimes\chi} M^* \otimes H\otimes M\otimes
H\xrightarrow{\id\otimes \alpha\otimes\id\otimes\beta \overline S}M\otimes M^*
\end{equation}
As before one easily verifies that $\mathrm{ev}^r$ and
$\mathrm{coev}^r$ are morphisms in $\bicomodH$ and also that they
define a right duality.

\begin{obse} \label{obse:dualsH} In explicit terms the comodule
  structures for the duals are given by the following formul\ae.
If $(M,\chi) \in \bicomodH$ and $f \in {}^*\!M$ and $m \in M$, then
${}^*\chi(f)=\sum f_{-1} \otimes f_{0} \otimes f_1 \in H \otimes {}^*\!M
\otimes H$ if and only if: $$\sum f_{0}(m) f_{-1} \otimes  f_1= \sum
f(m_0) \overline{S}(m_{-1}) \otimes S(m_1).$$ 
Similarly if $(M,\chi) \in \bicomodH$ and $f \in M^*$ and $m \in M$, then
$\chi^*(f)=\sum f_{-1} \otimes f_{0} \otimes f_1 \in H \otimes M^*
\otimes H$ if and only if: $$\sum f_{0}(m) f_{-1} \otimes  f_1= \sum
f(m_0) S(m_{-1}) \otimes \overline{S}(m_1).$$ 
\end{obse}

\begin{lema}
  For the category $\rcomodH$, the duality functors can be expressed in
  terms of the functors in Definition \ref{defi:()circrl}, in the
  following way.
$$
{}^*(-):(\rcomodH_f)^{\mathrm{op}}\xrightarrow{(-)^r}{}^H\!\M
\xrightarrow{(-)^\circ}\M^{H^{\mathrm{cop}}}_f
\xrightarrow{-\square_{H^{\mathrm{cop}}}S_+} \rcomodH_f
$$
$$
(-)^*:(\rcomodH_f)^{\mathrm{op}}
\xrightarrow{((-)^\circ)^{\mathrm{op}}}({}^{H^{\mathrm{cop}}}\!\M_f)^{\mathrm{op}}
\xrightarrow{(-)^\ell}\M^{H^{\mathrm{cop}}}_f
\xrightarrow{-\square_{H^{\mathrm{cop}}}\overline  S_+} \rcomodH_f
$$
\end{lema}
\begin{theo}\label{obs:dualsrlcirc2}
If $H$ is a finite dimensional coquasi Hopf algebra, then its
antipode has a canonical structure of a monoidal morphisms of coquasi
bialgebras $S:H^\circ\to H$. Moreover, this structure is given by the
functional $\chi^S:H\otimes H\to\Bbbk$ 
\begin{multline*}
\chi^S(x\otimes y)=\sum 
\phi^{-1}(S(y_3)\otimes S(x_3)\otimes
x_5)\alpha(x_4)\phi(S(y_2)S(x_2)\otimes x_6\otimes
y_5)\\
\alpha(y_4)\beta(x_8y_7)\phi(S(y_1)S(x_1)\otimes(x_7y_6)\otimes S(x_9y_8)).
\end{multline*}
and corresponds to the usual
monoidal structure of the left dual functor ${}^*(-)$. 
\end{theo}
\begin{proof}
By general categorical principles, the left dual functor has a
canonical monoidal structure
$^*(-):(\M^H)^{\mathrm{oprev}}\to\M^H$. This can be explicitly
computed in terms of the coquasi Hopf algebra structure of $H$. On the
other hand, we know the monoidal structures of the equivalences  $(-)^r$
and $(-)^\circ$, hence we can explicitly compute the monoidal
structure of $(-\square_{H^{\mathrm{cop}}}S_+)$. The latter is given
by a monoidal structure on the coquasi bialgebra morphism
$S:H^{\circ}\to H$ (see Theorem \ref{theo:fmonoidal}), 
and in fact it is given by the functional
$\chi^S$ above.
\end{proof}

Note
that the formula for $\chi^S$ above is written in function of the
comultiplication of $H$ not of the domain of $S$: $H^{\mathrm{cop}}$.

The functional in the theorem above appeared in \cite{kn:bulacu},
and in the dual case of quasi Hopf algebras in \cite{kn:dri}. 

\begin{theo}\label{theo:dualsS}
Let $H$ be a finite-dimensional coquasi Hopf algebra and consider the
monoidal structure on $S$ introduced in Proposition
\ref{obs:dualsrlcirc2} above. 
  The canonical linear isomorphisms $M\cong M^{\vee\vee}$ provide the
  components for monoidal natural isomorphisms 
$$
{}^{**}(-)\cong -\square_HS^2_+:\rcomodH_f\to\rcomodH_f \qquad
(-)^{**}\cong-\square_H \overline S^2_+:\rcomodH_f\to\rcomodH_f.
$$
\end{theo}
\begin{proof}
  If follows directly from Proposition \ref{obs:dualsrlcirc2}, Lemma
  \ref{l:relations} applied to $S:H^{\operatorname{cop}}\rightarrow H$ 
and Lemma \ref{l:roro}.
\end{proof}

\section{The fundamental theorem of Hopf modules}\label{s:fundtheo}
In this section we present a generalization to the case of
coquasi--Hopf algebras of the fundamental theorem
on Hopf modules introduced by Sweedler in \cite{kn:sweedler}. For a
modern presentation of this general case see also
\cite{kn:schau1,kn:Schauenburg}.  
\subsection{The fundamental theorem}
In this section we establish the basic set up in order to state and
prove the fundamental theorem on Hopf modules in our context. 

In this section, in order to use
(formal) duality arguments, we work with a braided monoidal category
that we call $\mathcal V$ and with unit object $\Bbbk$ (see
\cite{kn:JS} as a general reference). Moreover we assume that
$\mathcal V$ has equalizers and that the tensor product preserves 
equalizers in each variable.

The braiding $\gamma$ is a natural isomorphism
$\gamma_{X,Y}:X\otimes Y\to Y\otimes X$ and its existence ensures
that the natural transformations: 
\begin{multline*}
(X\otimes Y)\otimes (Z\otimes W )\xrightarrow{\cong}(X\otimes(Y\otimes
Z))\otimes W\xrightarrow{ (\id\otimes\gamma_{Y,Z})\otimes\id}
(X\otimes (Z\otimes Y))\otimes W\\\xrightarrow{\cong}(X\otimes Z)\otimes
(Y\otimes W)
\end{multline*}
define a monoidal structure on the functor 
$\otimes:\mathcal V\times\mathcal V\to \mathcal V$. Here we are
endowing the category $\mathcal V \times \mathcal V$ with its usual
monoidal structure: $(X,Y) \otimes (X',Y') = (X \otimes X', Y \otimes
Y')$.  

Using the coherence results in \cite{kn:JS}, 
we may assume without loss of generality
that $\mathcal V$ is a strict monoidal category. 

In the above set up, given a braided monoidal category $\mathcal V$, 
one can define coalgebra, comodule, coquasi bialgebra and
coquasi Hopf algebra objects in $\mathcal V$. The
braiding ensures that the tensor product of two coalgebras is a
coalgebra, and likewise with comodules. The fact that the tensor
product in $\mathcal V$ preserves equalizers in each variable, allows
us to define the cotensor product of bicomodules in exactly the same
manner than in the case that $\mathcal V$ is the category of
vector spaces. 

If $C$ is a coquasi bialgebra we  denote as 
$\mathcal V\,^C$, ${}^C\mathcal V $, ${}^C\mathcal V\,^C$,
${}^C_C\mathcal V\,^C$ and ${}^C\mathcal V\,^C_C$ 
the categories of right, left and bicomodules,
and the category of left and right Hopf modules in $\mathcal V$
respectively. The first three categories have monoidal
structures induced by the tensor product of $\mathcal V$. 

The category of bicomodules has also a
monoidal structure given by the cotensor product $\square_C$. 

In the same manner than in Definition \ref{defi:f+}, we
can define the bicomodules $f^+$ and $f_+$ 
for a morphism $f:C \rightarrow D$ of coalgebra objects in $\mathcal V$. 

For a coquasi bialgebra $C$ in $\mathcal V$, we call $u:\Bbbk
\rightarrow C$ the unit morphism.
In this situation we can consider a pair of adjoint functors 
$(u^+\square_\Bbbk-) \dashv (u_+\square_C
-):{}^C\mathcal V\,^C\to\mathcal V^C$. 

Explicitly, if $M$ is a bicomodule, $u_+\square_C
M$ is the right comodule of left coinvariants ${}^{\mathrm{co}C}M$ and if
$N$ is a right comodule, $u^+\square_\Bbbk N$ is the basic object 
$N \in  \mathcal V$ with the same
right coaction than $N$ and with 
left coaction given by $u\otimes\id_N:N\to C\otimes
N$. 
\begin{defi} \label{defi:free}
If $C$ is a coquasi bialgebra in $\mathcal V$, we 
define the free module functor $F:{}^C\mathcal V\,^C \rightarrow
{}^C_C\mathcal V\,^C$ as $F(M)=C \otimes M$ for $M \in {}^C\mathcal
V\,^C$. This functor, together with the forgetful functor
$U:{}^C_C\mathcal V\,^C\to {}^C\mathcal V\,^C$ constitute a pair of adjoint functors: $F\dashv U:{}_C^C\mathcal V\,^C
\to{}^C\mathcal V\,^C$.   
Define the functor $L:\mathcal V^C\to{}^C_C\mathcal V\,^C$ as the
composition $L= F \compo(u^+\square_\Bbbk-)$. Clearly $L$ will have a 
right adjoint given as $(u_+\square_C-)\compo  U$.
\end{defi}

The monoidal structure $\square_C$, lifts from ${}^C\mathcal V^C$
to the category of Hopf modules in such a way
that if $U: {}^C_C\mathcal V\,^C \rightarrow {}^C\mathcal V\,^C$ is the
forgetful functor, then the
adjunction $F\dashv U:{}^C_C\mathcal V\,^C\to {}^C\mathcal V\,^C$ is 
monoidal ({\em i.e.}, $U$ is lax monoidal, $F$ is strong monoidal and
the unit and counit of the adjunction are monoidal natural
transformations).

We want to show that under certain additional hypothesis, the functor
$L$ is a monoidal equivalence.

Without imposing further restrictions on the category $\mathcal V$, 
it is not hard to prove the following result, that is analogous
to \cite[Prop. 3.6]{kn:schau1} where the difference being that in the
above mentioned paper the result is formulated for the case that $\mathcal V=\mathbf{Vect}$ (see also  \cite[Lemma
2.1]{kn:schau1}).  
A general theorem that yields in particular the lemma we present
below, appears in \cite[Prop. 3.4]{kn:nacho}.

In this lemma we will need the following piece of notation.
If $C,D,C',D'$
are coalgebras in $\mathcal V$ and $U \in{}^C\,\mathcal V\,^D$ and 
$V\in{}^{C'}\,\mathcal V\,^{D'}$, we
will denote by $U\bullet V \in {}^{C\otimes C'}\,\mathcal V\,^{D\otimes
  D'}$ the object $U \otimes V$ equipped with the obvious bicomodule structure.

\begin{lema}\label{l:Lff}
  The functor $L: (\mathcal
  V\,^C,\Bbbk,\otimes) \rightarrow ({}^C_C\mathcal V\,^C,C,\square_C)$ 
is fully faithful and monoidal.
\end{lema}
 \begin{proof}
It is well known that the functor $L$ is fully faithful 
if and only if the unit of the adjunction
$L \dashv (u_+\square_C-)U$ is an isomorphism. It follows from the dual
of \cite[Lemma A1.1.1]{kn:johnstone}, that it is enough to exhibit a natural
isomorphism between $(u_+\square_C-)UL$ and the identity functor of
$\mathcal V^C$. 

The composition $UL:\mathcal V\,^C\to{}^C \mathcal V\,^C$ can be written
as $UL(M)=(C\bullet M)\square_{C^{\otimes 2}}p_+$. We have
natural isomorphisms 
$$
u_+\square_C (C\bullet M)\square_{C^{\otimes
    2}}p_+ \cong (u_+\bullet M)\square_{C^{\otimes
    2}}p_+ \cong M\square_C(u_+\bullet C)\square_{C^{\otimes 2}}p_+
\cong M
$$
where the last isomorphism is induced by $(u_+\bullet
C)\square_{C^{\otimes 2}}p_+\cong ((u\otimes \id_C)p)_+\cong (\id_C)_+=C$. 
This shows that there is a natural isomorphism $u_+\square_C UL(M)\cong M$. 

We will now exhibit a canonical monoidal structure on $L$.  The
basic observation is that, for $M\in\M^C$, $L(M)= C\otimes (u^+\square_\Bbbk M)$ is
isomorphic to $(C\bullet M)\square_{C^{\otimes 2}}p_+$, where $C\in
{}^C\!\M^C$ is the regular bicomodule. Then, we can form the
composition
\begin{align}
  L(M)\square_C L(N)&
  \cong ((C\bullet M)\square_{C^{\otimes
      2}}p_+)\square_C ((C\bullet M)\square_{C^{\otimes
      2}}p_+)\notag\\
  &\cong (((C\bullet M)\square_{C^{\otimes 2}}p_+)\bullet
  N)\square_{C^{\otimes 2}}p_+\notag\\
  &\cong(C\bullet M\bullet N)\square_{C^{\otimes 3}}(p_+\bullet
  C)\square_{C^{\otimes 2}}p_+ \label{eq:Lmonoidal1}\\
  &\cong(C\bullet M\bullet N)\square_{C^{\otimes 3}}(C\bullet p_+)
  \square_{C^{\otimes 2}}p_+\label{eq:Lmonoidal2}\\
  &\cong (C\bullet (M\otimes N))\square_{C^{\otimes 2}}p_+\notag\\
  &\cong L(M\otimes N)\notag.
\end{align}
All the isomorphisms above follow easily form the definition of the
contensor product except for the isomorphism between
\eqref{eq:Lmonoidal1} and \eqref{eq:Lmonoidal2}, which is induced by
the isomorphism $(p_+\bullet C)\square_{C^{\otimes 2}}p_+\cong
((p\otimes \id_C)p)_+\cong ((\id_C\otimes p)p)_+\cong (C\bullet
p_+)\square_{C^{\otimes 2}}p_+$ that is induced by the associator $\phi$. 

The isomorphism described above together with the obvious isomorphism
$\Bbbk\cong L(\Bbbk)$ provide a monoidal structure for $L$. The axioms of a
monoidal functor follow easily from the axioms satisfied by the 
associator $\phi$. 
\end{proof}

In the presence of an antipode, one obtains the following
strengthening of the above results. 
This form of the fundamental theorem on Hopf modules for
coquasi Hopf algebras is a consequence of \cite[Theorem
7.2]{kn:nacho}. It follows easily by a simple adaptation of the
arguments of  the Section 11
of the same work. 

\begin{theo}\label{t:fundtheo}
For an arbitrary  coquasi Hopf algebra in $\mathcal V$, the associated 
functor $L$ is a monoidal equivalence.   
\end{theo}
Observe that we do not ask $\mathcal V$ to be abelian or
additive. Neither we assume anything about the existence of duals in
$\mathcal V$. 
The only requirements on $\mathcal V$ are that it is braided monoidal,
with equalizers that are preserved by the tensor product.

For a version of Theorem \ref{t:fundtheo} over $\mathbf{Vect}$, see
for example \cite{kn:Schauenburg}.

In order to apply  this theorem to our context, we need its right
version that will be deduced below.

\begin{obse}\label{o:Vrev}
  Consider the braided monoidal category $\mathcal V^{\mathrm{rev}}$,
  which has the same underlying category as $\mathcal V$, the same
  unit object but the reverse tensor product
  $X\otimes^{\mathrm{rev}}Y=Y\otimes X$, see Section
  \ref{section:intro}. 
If $\gamma_{X,Y}: X \otimes Y \rightarrow Y \otimes X$ is the braiding
in $\mathcal V$, then the braiding in $\mathcal V^{\mathrm{rev}}$ is
$\gamma^{\mathrm{rev}}_{_{X,Y}}= \gamma_{_{Y,X}}$.  
The symmetry in the
  definition of coquasi bialgebra implies that if $(C,p,u,\phi)$ is a
  coquasi bialgebra in $\mathcal V$, then $(C,p,u,\phi)$ is a
  coquasi bialgebra in $\mathcal V^{\mathrm{rev}}$. Moreover, if
  $(S,\alpha,\beta)$ is an antipode for the coquasi bialgebra $H$ in
  $\mathcal V$, then $(S,\beta,\alpha)$ is an antipode for $H$ in
  $\mathcal V^{\mathrm{rev}}$. 
\end{obse}

\begin{coro}\label{c:fundth}
  Suppose $H$ is a coquasi Hopf algebra with invertible antipode in
  $\mathcal V$. Then the functor $R: {}^H \mathcal V
  \rightarrow {}^H\mathcal V^H_H$ defined as the composition of
  $ -\square_\Bbbk u_+ :{}^H\mathcal V\to{}^H\mathcal V^H$ with the
  free right Hopf module functor ${}^H\mathcal V^H\to {}^H\mathcal
  V^H_H$ is a monoidal equivalence from ${}^H\mathcal
  V$ to ${}^H\mathcal V^H_H$.
\end{coro}
\begin{proof}
  Let us denote $\mathcal V^{\mathrm{rev}}$ by $\mathcal W$. If
  $(H,u,p,\phi)$ is a coquasi bialgebra in $\mathcal V$, as we saw in
  Observation \ref{o:Vrev}, 
  $(H,u,p,\phi)$ is a coquasi bialgebra in $\mathcal W$, and if
  $(S,\alpha,\beta)$ is an antipode for $H$ in $\mathcal V$ then
  $(S,\beta,\alpha)$ is an antipode for $H$ in $\mathcal
  W$. 

  Clearly, as monoidal categories we have that 
  ${}^H\mathcal V=(\mathcal W^H)^{\mathrm{rev}}$,
  ${}^H\mathcal V^H=({}^H\mathcal W^H)^{\mathrm{rev}}$ and
  ${}^H\mathcal V^H_H=({}^H_H\mathcal W^H)^{\mathrm{rev}}$. 
The functor $-\square_\Bbbk u_+$ corresponds to
  $u^+\square-:\mathcal W^H\to{}^H\mathcal W^H $ and the free $H$-module
  functor ${}^H\mathcal V^H\to {}^H\mathcal V^H_H$ to the free
  $H$-module functor ${}^H\mathcal W^H\to {}^H_H\mathcal W^H$. Then
  $R$ is just $L^{\mathrm{rev}}$ for the coquasi Hopf algebra $H$ in
  $\mathcal W$, and hence it is a monoidal equivalence as we wanted to
  guarantee. 
\end{proof}

\subsection{The case of invertible antipode}
In this section we prove that the composition of the functors
$u^+\square_\Bbbk-:\M^H\to {}^H\!\M^H$ and the free right $H$--module
functor $F:{}^H\!\M^H\to {}^H\!\M^H_H$ is a monoidal equivalence when
$H$ has an invertible antipode. 

First we deal with the general case of a coquasi bialgebra.

\begin{lema} \label{lema:counit} For a coquasi bialgebra $C$, the functor 
$F(u^+\square_\Bbbk-): \rcomod  \rightarrow \hopfmod$ is left adjoint to 
$(u_+\square_H-)U: \hopfmod \rightarrow
\rcomod$ where $U:\hopfmod
\rightarrow \bicomod$ is the forgetful functor. 
Moreover, the counit transformation
corresponding to the above adjunction is the following: if $M
\in \hopfmod$, then $\varepsilon_M: {}_0({}^{\mathrm{co}C}M) \otimes C
\rightarrow M$ is the map $\varepsilon_M(m \otimes c)= m \cdot c$.
\end{lema}
\begin{proof}
The assertion about the adjunction is clear. If $\varepsilon'$ and
$\varepsilon''$ are the counits of $F\dashv U$ and
$(u^+\square_\Bbbk-)\dashv (u_+\square_C-)$ respectively, then the
counit of the composition of these functors is 
$$
\varepsilon:F (u^+\square_\Bbbk-)(u_+\square_C-)U\xrightarrow{F\varepsilon''
  U}FU\xrightarrow{\varepsilon '}\id.
$$
For $M\in {}^C\!\M^C$, the transformation $\varepsilon''_M$ is the
inclusion of ${}_0({}^{\mathrm{co}C}M)$ in $M$, while for $N\in
\hopfmod$, $\varepsilon '_N:N\otimes C\to N$ is given by the right action
of $C$ on $N$. Therefore $\varepsilon $ is indeed given 
by the above formula. 
\end{proof}

\begin{defi}\label{defi:iota}
Let $H$ be a coquasi Hopf algebra. 
Define a functor  ${\mathcal{I}}:\M^H\to{}^H\!\M$ as the composition
$\rcomodH\xrightarrow{(-)^\circ}{}^{H^{\mathrm{cop}}}\!\M\xrightarrow{S^+
\square_{H^{\mathrm{cop}}}-}\lcomodH$. 
In other words, on objects  ${\mathcal{I}}(M,\chi)=(M,(S \otimes
\id)\mathrm{sw}\chi)$, and on arrows  ${\mathcal{I}}$ is the identity.
\end{defi}

\begin{coro}\label{coro:iotamonoidal}  In the situation above the
  functor ${\mathcal{I}}:\mathcal M^H{^{\mathrm{rev}}} \rightarrow \lcomodH$
  is monoidal 
\end{coro} 
\begin{proof} It follows immediately from:\,a) the equality 
\[{\mathcal{I}} = (S^+ \square_{H^{\mathrm{cop}}}-)\compo  (-)^\circ :
\rcomodH\xrightarrow{(-)^\circ}{}^{H^{\mathrm{cop}}}\!\M\xrightarrow{S^+
  \square_{H^{\mathrm{cop}}}-}\lcomodH ; \] b) the fact that $(-)^\circ$ is
monoidal --see the comments after Definition \ref{defi:()circrl}--; 
c) Theorem \ref{theo:fmonoidal}.  
\end{proof}

The following natural transformation will be crucial in the proof of
Radford's formula.

\begin{theo} \label{lema:tau} 
  Let $H$ be a coquasi Hopf algebra with invertible antipode. 
For $M \in \rcomodH$ the arrows
  $\tau_M: M \otimes H \rightarrow M \otimes H$ defined as $\tau_M(m
  \otimes h) = \sum m_0 \phi^{-1}(m_1 \otimes S(m_3) \otimes h_2)\beta(m_2)
  \otimes S(m_{4})h_1$ are the components of a natural 
transformation between the functors  $F \compo\, (-)_0\, \compo\, {\mathcal{I}}$
and 
$F \compo\ _0(-):\rcomodH \rightarrow \hopfmodH$. Moreover, the natural
transformation $\tau$ is 
invertible and its inverse is given for all $M \in \rcomodH$ by the
formula: $\tau^{-1}_M(m \otimes h)=\sum
\phi(S(m_1)\otimes m_3\otimes h_1)\alpha(m_2)m_0\otimes m_4h_2$.
\end{theo}
\begin{proof}

It is convenient for the proof to split the map $\tau_M$ as
follows. First define the map $\pi_M: ({\mathcal{I}} M)_0 \rightarrow {}_0M
\otimes H$ as $\pi_M(m)=\sum m_0 \otimes \beta(m_1)S(m_2)$. An
elementary computation shows that 
\begin{equation*}\tau_M:({\mathcal{I}} M)_0 \otimes H
\xrightarrow{\pi_M \otimes \id} ({}_0M \otimes H) \otimes H 
\xrightarrow{\Phi_{{}_0M,H,H}} 
{}_0M \otimes (H \otimes H)\xrightarrow{\id \otimes p} {}_0M \otimes H
\end{equation*}

 The structure of $H$--bicomodule on $(F \compo\,(-)_0
  \,\compo\, {\mathcal{I}})(M)=M\otimes H$ is $m \otimes h \mapsto \sum S(m_1)h_1
  \otimes m_0 \otimes h_2 \otimes h_3$, while
  the structure on $(F \compo\ _0(-))(M)=M \otimes H$ is $m \otimes h
  \mapsto \sum h_1 \otimes m_0 \otimes h_2 \otimes m_1h_3$. A
  direct verification shows that $\pi_M$ is a morphism of
  bicomodules. Hence $\tau_M$ is a composition of morphisms of
  bicomodules. The compatibility of $\tau_M$ with the right action of $H$ is
  deduced directly from the fact that $\pi_M \otimes \id$ and 
$(\id \otimes p)\Phi_{{}_0M,H,H}$ are morphisms of right
$H$--modules. The $H$--equivariance of the first morphism is obvious,
while the equivariance of the second is a consequence of the following
general fact --that we apply in the situation that  
$\mathcal C=\bicomodH$ and $A=H$--. 
If $(A,p,u)$ is an algebra --also called a {monoid}--
in an arbitrary monoidal category $\mathcal C$ (that in accordance with
\cite{kn:JS} can be assumed to be strict) then, for any object $X$
the arrow $\id\otimes p:X\otimes A\otimes A\to X\otimes A$ is a
morphism of right $A$--modules. 

Finally, the verification of that the maps $\tau_M$ and $\tau_M^{-1}$
are indeed inverses to each other is a direct computation. 
\end{proof}
 
A version of the above lemma for quasi Hopf algebras appears in
\cite{kn:Schauenburg}. Our proof is similar. 

The following result is an 
immediate consequence of Theorem \ref{lema:tau} and Corollary
\ref{c:fundth}. 
\begin{coro}\label{c:ft}
  In the situation above, the functor $F(u^+\square_\Bbbk -):(\mathcal
  M^H)^{\mathrm{rev}}\to{}^H\mathcal M^H_H$ has a unique monoidal
  structure such that $\tau$ is a monoidal natural transformation. In
  particular, with this structure $F(u^+\square_{\Bbbk}-)$ is a
  monoidal equivalence. 
\end{coro}

We end the section with the following observation, that will be used
in Section \ref{ss:monoidality}. 
\begin{obse}\label{o:fsquareg}
  If $M,N\in\M^H$ and $P,Q\in {}^H\!\M$, there exist canonical
  isomorphisms of bicomodules
  $$
  ((u^+\square_\Bbbk M)\otimes H)\square_H ((u^+\square_\Bbbk N)\otimes H)
  \cong (u^+\square_\Bbbk N)\otimes((u^+\square_\Bbbk M)\otimes H)
  $$
  $$
  ((P\square_\Bbbk u_+)\otimes H)\square_H((Q\square_\Bbbk u_+)\otimes H)
  \cong (P\square_\Bbbk u_+)\otimes ((Q\square_{\Bbbk}u_+)\otimes H)
  $$

If $f:(P\square_\Bbbk u_+)\to (u^+\square_\Bbbk M)$ and $g:
(Q\square_\Bbbk u_+)\to (u^+\square_\Bbbk N)$ are morphisms of
bicomodules, then the following diagram commutes
$$
\xymatrixrowsep{.7cm}
\diagramcompileto{pq}
  ((P\square_\Bbbk u_+)\otimes H)\square_H((Q\square_\Bbbk u_+)\otimes
  H)
\ar[r]^-\cong \ar[dddd]_{f\square_Hg}
&
(P\square_\Bbbk u_+)\otimes((Q\square_\Bbbk u_+)\otimes H)
\ar[d]^{\id\otimes g}
\\
&(P\square_\Bbbk u_+)\otimes((u^+\square_\Bbbk N)\otimes H)
\ar[d]^{\Phi^{-1}}
\\
&((P\square_\Bbbk u_+)\otimes(u^+\square_\Bbbk N))\otimes H
\ar[d]^{\mathrm{sw}\otimes \id}
\\
&((u^+\square_\Bbbk N)\otimes (P\square_\Bbbk u_+))\otimes H
\ar[d]^{\id\otimes f}
\\
((u^+\square_\Bbbk M)\otimes H)\square_H ((u^+\square_\Bbbk N)\otimes
H)
\ar[r]^-\cong
&
(u^+\square_\Bbbk N)\otimes ((u^+\square_\Bbbk M)\otimes H)
\enddiagram
$$
\end{obse}

\section{The Frobenius isomorphism and the object of cointegrals}
\label{s:frobiso}
Suppose that $H$ is a coquasi Hopf algebra with invertible antipode.
If $M$ is a left $H$--module in the category $\bicomodH_f$
  its left dual ${}^*M$, is an object in $\hopfmodH$ in a functorial
  way.  If $a_M: H \otimes M \rightarrow M$ is a $H$--module
  structure of $M$ in $\bicomodH$, the corresponding structure
  $a_{{}^*M}: {}^{*}M \otimes H \rightarrow {}^*M$ is given by
\begin{equation}\label{e:dualstructure}
\begin{split}
{}^*M\otimes H\xrightarrow{\id\otimes\mathrm{coev}^\ell}&({}^*M\otimes
H)\otimes (M\otimes {}^*M)\xrightarrow{\cong}({}^*M\otimes( H\otimes
M))\otimes  {}^*M\\ 
\xrightarrow{(\id\otimes a_M)\otimes\id}&({}^*M\otimes M)\otimes
{}^*M\xrightarrow{\mathrm{ev}^\ell\otimes 1}{}^*M. 
\end{split}
\end{equation}

From now on we assume that $H$ is a finite dimensional coquasi Hopf algebra. 
If we take $H \in \bicomodH_f$ as a left
$H$--module 
with respect to the regular action, its right dual ${}^*H$ is canonically an
object in $\hopfmodH$. An explicit description of the right 
$H$--structure defined above for ${}^*H$ is the following: if $f\in{}^*H$
and $x,y\in H$, 
$a_{{}^*H}(f\otimes x)(y)=(f \cdot x)(y)$ is equal to

\begin{equation}\label{e:fourier}
\begin{split}
(f \cdot x)&(y)=\\
\sum & \phi^{-1}(\overline{S}(x_5y_7)x_1\otimes y_3\otimes \overline
S(y_1))\alpha\overline{S}(y_2)\phi(\overline{S}(x_4y_6)\otimes x_2\otimes
y_4)\beta\overline{S}(x_3y_5)\\ f(x_6&y_8) \phi(S(x_7y_9)x_{11} \otimes
y_{13} \otimes S(y_{15}))\phi^{-1}(S(x_8y_{10})\otimes x_{10} \otimes
y_{12})
\alpha(x_9y_{11})\beta(y_{14}).
\end{split}
\end{equation}

It is important to notice that in the above formula --and in the
formula for the Frobenius isomorphism-- we obtain the expression for
$f \cdot x \in {}^*H$ in
terms of the {\em standard } evaluation of vector spaces
$H^{\vee}\otimes H\to \Bbbk$.  

\begin{theo}\label{lema:basic}
If $H$ is a finite dimensional coquasi Hopf algebra, 
then there exists a unique up to
  isomorphism one dimensional object $W \in \rcomodH$ such that 
there is an isomorphism ${}_0W \otimes H \cong {}^*H \in \hopfmodH$. 
Moreover, $W$ can be taken
as the space of left cointegrals of the Hopf algebra $H$ and the
isomorphism --called the Frobenius isomorphism-- is the map 
$\mathcal F$ given by 
\begin{equation}\label{e:fouriertransform0}
\mathcal F(\varphi \otimes x) = \varphi \cdot x
\end{equation}
where the action used is the one defined in \eqref{e:dualstructure} 
and applied to $M={}^*H$ in accordance to the formula \eqref{e:fourier}. 
\end{theo}
\begin{proof} The existence and uniqueness of $W$ follows
  immediately from the fundamental theorem on Hopf modules--see more
  specifically Corollary \ref{c:fundth}--. The
  characterization of $W$ as a space of left cointegrals  is deduced directly
  from the explicit description of the inverse functor of $F\compo
  {}_0(-)$ as the composition of the 
forgetful functor $U: {}^H\M^H_H \rightarrow \M^H_H$ with the left
fixed part functor --see the considerations previous to Lemma
\ref{l:Lff}--. Thus,  
  $W$ is the space of left coinvariants of\,
  $\/{}^*H$ with respect to the coaction described in
  \eqref{e:*Mstructure}. In explicit terms $W=\{\varphi \in {}^*H:
  {}^*\chi(\varphi)= 1 \otimes \varphi\}$. 
Using the description of ${}^*\chi$, appearing
  in Observation \ref{obse:dualsH}, we conclude that $\varphi \in W$ if and
  only if for all $x \in H$, $\varphi(x)1 = \sum \overline{S}(x_1)
  \varphi(x_2)$.  In other words $\varphi(x)1 = \sum x_1 \varphi(x_2)$
and then $\varphi
  \in {}^*H$ is a left cointegral.  
 The description of the
  counit of the adjunction as the map given by the action --see Lemma \ref{lema:counit}-- will
    yield the characterization \eqref{e:fouriertransform0}.
\end{proof}

\begin{obse} \label{o:modular} In the same manner than in the
  classical case, from the existence of the
  isomorphism $\mathcal F$ we
conclude that $W$, the space of left cointegrals,  is one dimensional.
Hence, one can prove 
the existence of a group like element $a\in H$ such that
$\sum \varphi(x_1)x_2=\varphi(x)a$ for all $\varphi\in W$ and $x\in
H$. The element  $a \in H$ is called the \emph{modular} element.
\end{obse}
\begin{obse}\label{o:explicit}
Using the definition of the modular element $a$ just presented as well
as formula \eqref{e:fourier} applied to the situation that $f =
\varphi$ is a cointegral, we obtain the following explicit formula for
the Frobenius isomorphism.
\begin{equation}\label{e:fouriertransform}
\begin{split}
\mathcal F(\varphi & \otimes x)(y)=\alpha(a)\beta(1)\\ \sum & \phi^{-1}(x_1
\otimes y_2 \leftharpoonup \alpha \otimes \overline{S}y_1) \varphi(x_2y_3)
\phi(x_3 \otimes \beta \rightharpoonup y_4 \otimes
Sy_5) \phi(a^{-1} \otimes a \otimes Sy_6)
\end{split}
\end{equation}

\end{obse}
\begin{lema}
  The coaction $\chi_W:W\to W\otimes H$ is of the form
  $\chi_W(\varphi) = \varphi \otimes a^{-1}$ for $a \in H$ as above. 
\end{lema}
\begin{proof}
As $W$ is one dimensional, the coaction  $\chi_W$ is of the form
$\chi_W(\varphi) = \varphi \otimes b$ for $b \in H$. 
The definition of the right
comodule  structure on $W$ (see Observation \ref{obse:dualsH})
yields --for $x \in H$-- the formula: $\sum
\varphi(x_1)S(x_2)=\varphi(x)b$. 
It follows then that $S(a)=a^{-1}=b$ --see Observation
\ref{obse:grouplike}--.  
\end{proof}
\begin{obse}\label{W=b*}
\noindent a) The comodule $W$ is isomorphic to ${a^{-1}}_+\in\rcomodH$, where 
$a^{-1}:\Bbbk\to H$ is the coalgebra morphism induced by the 
multiplication by $a^{-1}$.  

\noindent b) Similarly, the coaction in ${}^{*}W$ is given as 
$\chi_{_{{}^{*}W}}(t) = t \otimes a$ for  $t \in {}^{*}W$. Hence,
${}^{*}W$ is isomorphic to $a_+\in\rcomodH$. 
\end{obse}

Recall that we abbreviated the 
functors $u^+\square_\Bbbk-$ and $-\square_{\Bbbk}u_+$  by ${}_0(-)$
and $(-)_0$ respectively. 

\begin{obse}\label{obse:W}
  $W \in
  \rcomodH$ as well as ${}_0W \in \bicomodH$ are invertible objects in the
  corresponding monoidal categories. In other words, the functors $-
  \otimes W: \rcomodH \rightarrow \rcomodH$ and 
   $- \otimes {}_0W: \bicomodH \rightarrow \bicomodH$ are equivalences
   and it is clear that  
  the inverse equivalences are obtained by
  tensoring with the corresponding duals.    
\end{obse}

For use in the next section we write down the following definitions.
\begin{defi}\label{defi:c_W} Define the following functors
  $c_W^{l},c_W^{r}:\bicomodH \rightarrow \bicomodH$ as follows: $c_W^{l}=
  ({}_0W \otimes -) \otimes {}_0{}^{*}W$ and  $c_W^{r}=
  {}_0W \otimes (- \otimes {}_0{}^{*}W)$  
\end{defi}

\begin{obse}\label{obse:cmonoidal} It is clear that $c_W^{l}$ and
  $c_W^{r}$ are monoidal functors that are naturally isomorphic
  via the natural transformation given by the obvious associator. 
In the notations of Theorem \ref{lema:basic} and using the fact that
$c_W^{l},c_W^{r}$ are monoidal functors, we conclude that 
$({}_0W \otimes
H) \otimes {}^*_0W$  and $({}_0W \otimes H) \otimes {}^*_0W$ are
algebras in the category $\bicomodH$. 

\end{obse}
\section{Radford's formula}\label{s:formula}

In this section we use categorical methods to prove Radford's formula 
expressing $S^4$ in terms of conjugation with a functional and a group
like element. In the second part of this section we prove the
monoidality of the functional.

\subsection{Radford's formula}
We use the notations of the last section and assume that $H$ is a
finite dimensional coquasi Hopf algebra.
We will take basis elements $\varphi \in W$, $t \in {}^*W$ normalized
in such a way that $t(\varphi)=1$. 

\begin{lema}\label{lema:gamma} In the notations of Theorem \ref{lema:basic} the
  isomorphism in $\bicomodH$ 
$$\gamma: ({}_0W \otimes H) \otimes {}^*_0W
  \xrightarrow{\mathcal F \otimes \operatorname {id}}{}^*H \otimes {}^*_0W \cong
  {}^*({}_0W \otimes H)\xrightarrow{({}^*\mathcal{F})^{-1}}
  {}^{**}H$$ 
is a morphism of algebras. Moreover, if we define the Nakayama isomorphism
$\mathcal N:H \rightarrow {}^{**}H$ by the formula: $\mathcal
N(x)=\gamma((\varphi \otimes x) \otimes t)$, then the commutativity of the
diagram below characterizes $\mathcal N$:
\[
\diagramcompileto{lemma:gamma0}
H \otimes (W \otimes H)\ar[d]_{\mathcal N \otimes
  \mathcal F} \ar[rr]^{\cong}&& 
(H \otimes W) \otimes H \ar[d]^{\mathcal F \operatorname{sw} \otimes
  \operatorname {id}}\\
 {}^{**}H \otimes {}^*H\ar[r]_-{\operatorname{ev}_{{}^*H}^{\ell}} &\Bbbk& {}^*H \otimes H\ar[l]^-{\operatorname{ev}_{H}^{\ell}}
\enddiagram
\]
\end{lema}
\begin{proof} 
The multiplicativity of $\gamma$ follows immediately from the
fact that $\mathcal F$ is a morphism of $H$--modules and from the
commutativity of the following diagram that is a direct consequence of the
definition of the action $a_{{}^*H}: {}^*H \otimes H \rightarrow
{}^*H$ --see Definition \eqref{e:dualstructure}--.
\[
\diagramcompileto{lemma:gamma}
({}^*H \otimes H) \otimes H\ar[d]_{a_{{}^*H}\otimes 
  \operatorname {id}} \ar[rr]^{\cong}&& 
{}^*H \otimes (H \otimes H) \ar[d]^{
  \operatorname {id}\otimes p}\\
 {}^{*}H \otimes H\ar[r]_-{\operatorname{ev}_{H}^{\ell}} &\Bbbk& {}^*H \otimes H\ar[l]^-{\operatorname{ev}_{H}^{\ell}}
\enddiagram
\]
The assertion
concerning $\mathcal N$ follows directly from the definitions. 
\end{proof}
\begin{obse}\label{obse:generalgamma} 
\noindent
a) The fact that $\gamma$ is a
  morphism of algebras is valid in the following general context. Let
  $\mathcal C$ be a rigid monoidal category and let $a,w \in \mathcal C$
  be respectively an algebra  and an arbitrary object.
 Let $\mathcal F: w \otimes a \rightarrow {}^*a$ be an invertible 
morphism of $a$--modules in $\mathcal C$, then the 
object $(w \otimes a) \otimes {}^*w$ is an algebra in $\mathcal C$ and
the map 
$$\gamma: (w \otimes a) \otimes {}^*w
  \xrightarrow{\mathcal F \otimes \operatorname {id}}{}^*a \otimes
  {}^*w \cong {}^*(w \otimes a)\xrightarrow{({}^*\mathcal{F})^{-1}}
  {}^{**}a$$ is a morphism of algebras. 

\noindent
b) In the case of ordinary Hopf algebras, the commutativity of the diagram that
characterizes $\mathcal N$ after 
identifying $H$ with its double dual reads as $\varphi(y{\mathcal
  N}(x))= \varphi(xy)$ that is the usual definition of the 
Nakayama automorphism.
\end{obse}     
\begin{lema} \label{lema:xi} In the notations of Theorem \ref{lema:tau} 
and Theorem  \ref{lema:basic} if $M$ is an
    object in $\rcomodH$, then the morphism  $\xi_M: {}^{**}{\mathcal{I}}(M)_0\otimes 
(({}_0W\otimes H)\otimes
{}^*_0W)\to  {}^{**}{}_{\hspace{.1cm}0}M \otimes 
(({}_0W\otimes H)\otimes
{}^*_0W)$ defined
by the commutativity of the diagram below is a
morphism of right $({}_0W\otimes H) \otimes {}^*_0W$--modules. 
\begin{equation}\label{e:xi}
\begin{split}
\diagramcompileto{defxidiag}
{}^{**}({\mathcal{I}}(M)_0\otimes H)\ar[d]_\cong
\ar@{->}[rr]^-{{}^{**}\tau_M}&&{}^{**}({}_0M\otimes H)\ar[d]^\cong \\
{}^{**}{\mathcal{I}}(M)_0\otimes
{}^{**}H\ar[d]_{\operatorname{id} \otimes\gamma^{-1}}\ar@{->}[rr]^{\omega_M} && 
{}^{**}_{\hspace{.1cm}0}M\otimes
{}^{**}H\ar[d]^{\operatorname{id}\otimes\gamma^{-1}}
\\ {}^{**}{\mathcal{I}}(M)_0\otimes 
(({}_0W\otimes H)\otimes
{}^*_0W) \ar@{->}[rr]^-{\xi_M}&&
 {}^{**}_{\hspace{.1cm}0}M\otimes 
(({}_0W\otimes H)\otimes
{}^*_0W)
\enddiagram
\end{split}
\end{equation}
\end{lema}
\begin{proof}
Being $\tau_M$ a morphism of $H$--modules it is
  clear that $\omega_M$ is a morphism of ${}^{**}H$--modules. Then,
  from the fact that $\gamma$ is an algebra morphism and that all the
  modules involved are free over the corresponding algebra objects, it
  follows that $\xi_M$ is a morphism of $(W_0 \otimes H) \otimes
  {}^{*}W_0$--modules.
\end{proof}

\begin{defi}\label{defi:nueme}
  Define $\nu_M: {}^{**}{\mathcal{I}}(M)_0\otimes ({}_0W\otimes H)\to{}^{**}
  _{\hspace{.1cm}0}M\otimes ({}_0W\otimes H)$ as the unique morphism such
  that the diagram below commutes. 
\[
\diagramcompileto{nudef}
{}^{**}{\mathcal{I}}(M)_0\otimes 
(({}_0W\otimes H)\otimes
{}^*_0W) \ar@{->}[rr]^-{\xi_M}\ar[d]_\cong&&
 {}^{**}_{\hspace{.1cm}0}M\otimes 
(({}_0W\otimes H)\otimes
{}^*_0W)\ar[d]^\cong\\
({}^{**}{\mathcal{I}}(M)_0\otimes 
({}_0W\otimes H))\otimes
{}^*_0W\ar@{->}[rr]^-{\nu_M\otimes \operatorname{id}_{{}^*_0W}}&&
({}^{**}_{\hspace{.1cm}0}M\otimes 
({}_0W\otimes H))\otimes
{}^*_0W
\enddiagram
\]

  The existence of
  $\nu_M$ and the fact that it is a
  morphism in $\bicomodH$ follows immediately from
  the considerations of Observation \ref{obse:W}. Moreover from the
  fact that  $\xi_M$ is a morphism in
  ${}^H\mathcal M^H_{({}_0W\otimes H)\otimes {}_0^*W}$, it follows
  that $\nu_M$
    is a morphism in   $\hopfmodH$. 

The monoidality of ${\mathcal{I}}$ gives canonical isomorphisms
${\mathcal{I}}(M^{**})\cong{}^{**}{\mathcal{I}}(M)$. Composing with
$\nu_M$ we get an isomorphism 
$\widehat\nu_M:{\mathcal{I}}(M^{**})_0\otimes ({}_0W\otimes H)\to{}^{**}
  _{\hspace{.1cm}0}M\otimes ({}_0W\otimes H)$.
\end{defi}

For $M\in\rcomodH$, consider the following composition or arrows in
$\hopfmodH$, from ${}_0W^*\otimes (({}^{**}_0\! M \otimes {}_0 W)\otimes
H)$ to ${}_0M^{**}\otimes H$.
\begin{multline}
  \label{eq:zeta2}
  \zeta_M:{}_0W^*\otimes {}^{**}_0\! M \otimes {}_0 W\otimes H
  \xrightarrow{\id\otimes\widehat\nu^{-1}_M} {}_0W^*\otimes
  {\mathcal{I}}(M^{**})_0\otimes {}_0W\otimes H\to\\
  \xrightarrow{\id\otimes\mathrm{sw}\otimes\id} {}_0W^*\otimes
  {}_0W\otimes {\mathcal{I}}(M^{**})_0\otimes
  H\xrightarrow{\mathrm{ev}\otimes\tau_{M^{**}}} {}_0M^{**}\otimes H
\end{multline}
Here we omitted the associativity constraints for simplicity. However
this does not introduce any ambiguity as long as we know how to
associate the domain and codomain, by the coherence theorem for
monoidal categories. 

The composition $\zeta_M$ is a morphism in $\hopfmodH$. Indeed, $\widehat\nu_M$
and $\tau_{M^{**}}$ are morphisms of Hopf modules; the morphism
  $\id\otimes \mathrm{sw}\otimes \operatorname{id}$ is the image under
  the free Hopf module functor 
  $\bicomodH\to\hopfmodH$ of the morphism of bicomodules
  $\id\otimes \mathrm{sw}:{}_0W^{*}\otimes {\mathcal{I}}(M^{**})_0\otimes
  {}_0W\to {}_0W^*\otimes {}_0 W\otimes 
  {\mathcal{I}}(M^{**})_0$. Observe that $\mathrm{sw}$ is a morphism
  of bicomodules 
  because the trivial comodule structures in each tensor factor are
  added on opposite sides. 

\begin{defi} \label{defi:zeta}
    Denote by $\mu_M:W^*\otimes ({}^{**} M\otimes W)\to M^{**}$ the unique
    morphism in $\rcomodH$ such that $\mu_M\otimes \id_{H}=\zeta_M$.
\end{defi}

It is clear that $\mu$ is a natural isomorphism between the functors 
$W^*\otimes( {}^{**}(-)\otimes W)$ 
and $(-)^{**}: \rcomodH \rightarrow \rcomodH$. 

\begin{coro}\label{c:Msquare()+}
  The canonical linear isomorphisms $M\cong M^{\vee\vee}$ together
  with $\mu$ give a natural isomorphism in $\rcomodH$ 
\begin{equation}\label{e:Msquare()+}
M\square_H\, p(a^{-1}\otimes p(\overline S^2\otimes a))_+\to M\square_H S^2_+.
\end{equation}
\end{coro}
\begin{proof}
First we use Theorem \ref{theo:dualsS}, and substitute ${}^{**}M$ by
$M \square_HS^2_+$ and $M^{**}$ by $M \square_H \overline S^2_+$. In
this manner we obtain from $\mu_M$ an isomorphism
$$
W^*\otimes ((M\square_H\overline S^2_+)\otimes W)\to M\square_HS^2_+.
$$
Now using the fact that $W\cong a^{-1}_+$ --Observation \ref{W=b*}--
and the conclusions of Observation \ref{obse:squareandtensor} part b) we deduce our
result.
\end{proof}

\begin{theo}[Radford's formula]\label{theo:radford}
  There exists an invertible functional $\sigma:H\to\Bbbk$ such that
  for all $x\in H$ 
$$
a^{-1}(\overline S^2(x)a) =S^2(\sigma \rightharpoonup x \leftharpoonup
\sigma^{-1}). 
$$
\end{theo}
\begin{proof}
  It easily follows form Corollary \ref{c:Msquare()+} and Theorem
  \ref{l:nacho}. Indeed in the  situation of Corollary
  \ref{c:Msquare()+} the theorem guarantees the
  existence of a functional $\sigma$ such that --see Observation
  \ref{obse:concrete}, \eqref{e:concrete}--
  $p(a^{-1}\otimes p(\overline S^2\otimes a))(x)=S^2
(\sigma \rightharpoonup x \leftharpoonup\sigma^{-1})$.
\end{proof}

The functional $\sigma$ defined in the theorem above is the analogue
for finite dimensional coquasi Hopf algebras of the modular function
of a finite dimensional Hopf algebra. See Section \ref{s:casehopf}.

\begin{obse}\label{obse:final} The above formula can be transformed
  into another similar to the classical formula:
\begin{equation}\label{e:final4} 
S^4(x) = (a^{-1}(\widehat{\sigma}\rightharpoonup x \leftharpoonup
\widehat{\sigma}^{-1}))a 
\end{equation}
where $\widehat{\sigma}$ is another invertible functional that can be
computed explicitly in terms of the above information.
\end{obse}

\subsection{Monoidality}\label{ss:monoidality}

In this section we prove that the natural isomorphism $\mu$ of
Definition \ref{defi:zeta} is monoidal. We shall work as if the
monoidal category $({}^H\!\M^H,\Bbbk,\otimes)$ were strict, and hence
ignore the associativity and unit constraints. This can be formalized
by passing to an monoidally equivalent strict monoidal
category. Indeed, 
our proof does not depend on the fact that we are working with
the category of comodules, but only on certain properties satisfied by
the several arrows we consider. 

The functor $\M^H\to \M^H$ given by $M\mapsto W^*\otimes
M\otimes W$ has a canonical monoidal structure given by the
constraints
$$
W^*\otimes M\otimes W\otimes W^*\otimes N\otimes
W\xrightarrow{\id\otimes\id\otimes\mathrm{ev}\otimes\id\otimes\id}
W^*\otimes M\otimes N\otimes W
$$
$$
\Bbbk\xrightarrow{\mathrm{coev}}W^*\otimes
W\xrightarrow{\cong}W^*\otimes\Bbbk\otimes W
$$
These morphisms are isomorphisms because $W$ is an invertible object
--it has dimension one--.  

\begin{theo}\label{t:mumon}
  The natural transformation $\mu$ in Definition \ref{defi:zeta} is monoidal.
\end{theo}
The assertion that $\mu$ is a monoidal natural transformation is
expressed in the commutativity of the diagrams in Figure
\ref{fig:mumon0}.
\begin{figure}
  \begin{equation}\label{eq:mumon0.a}
  {\xymatrixcolsep{1.2cm}
  {\diagramcompileto{mumon0}
  W^*\otimes {}^{**}M\otimes W\otimes W^*\otimes {}^{**}N\otimes W
  \ar[r]^-{\mu_M\otimes\mu_N}
  \ar[d]|-{\id\otimes\id\otimes\mathrm{ev}\otimes\id\otimes\id\hole\hole\hole}
  &
  {}^{**}M\otimes {}^{**}N
  \\
  W^*\otimes M\otimes N\otimes W
  \ar[ur]_{\mu_{M\otimes N}}
  &
  \enddiagram}}
\diagramcompileto{eq:mumon0.b}
&\Bbbk\ar[d]^\id\ar[dl]\\
W^*\otimes\Bbbk\otimes W\ar[r]_-{\mu_\Bbbk}
&
\Bbbk
\enddiagram
  \end{equation}
  \caption{ }
  \label{fig:mumon0}
\end{figure}

\begin{proof}
We divide the proof in two parts. In some diagrams, we omit the symbol
$\otimes$ as a saving space measure, adding parenthesis when
necessary. 
\begin{figure}
$$ 
  \def\objectstyle{\scriptscriptstyle}
  \def\labelstyle{\scriptscriptstyle}
  \xymatrixcolsep{.4cm}
  {\diagramcompileto{mumon1}
  ({}_0W^*)({}_{\phantom{*}0}^{**}M) ({}_0W)
  ({}_0W^*)({}_{\phantom{*}0}^{**}N)({}_0W) H
  \ar[r]^-{(\id)\mathrm{ev}(\id)}
  \ar[d]_{(\id)\hat\nu_N^{-1}}
  &
  ({}_0W^*)({}_{\phantom{*}0}^{**}M)({}_{\phantom{*}0}^{**}N)({}_0W)
  H
  \ar[r]^-\cong
  \ar[d]_{(\id)\hat\nu_{N}^{-1}}
  \ar@{}[ddddr]|-{\scriptstyle\mathrm{(C)}}&
  ({}_0W^*)({}_{\phantom{*}0}^{**}(MN))({}_0W) H
  \ar[dd]^{(\id)\hat\nu_{MN}^{-1}}
  \\
  ({}_0W^*)({}_{\phantom{*}0}^{**}M) ({}_0W)
  ({}_0W^*)({\mathcal{I}}(N^{**}))({}_0W) H
  \ar[r]\ar[d]_{(\id)\mathrm{sw}(\id)}
  &
  ({}_0W^*)({}_{\phantom{*}0}^{**}M)({\mathcal{I}}(N^{**}))({}_0W) H 
  \ar[d]|-{(\id)\mathrm{sw}(\id)}
  &
  \\
  ({}_0W^*)({}_{\phantom{*}0}^{**}M) ({}_0W)
  ({}_0W^*)({}_0W)({\mathcal{I}}(N^{**}))H 
  \ar[d]_{(\id)(\id)(\id)\mathrm{ev}(\id)(\id)}\ar[r]
  &
  ({}_0W^*)({}_{\phantom{*}0}^{**}M)({}_0W)({\mathcal{I}}(N^{**})_0) H
  \ar[d]|-{\mathrm{sw}(\id)}
  &
  ({}_0W^*)({\mathcal{I}}(MN)_0)({}_0W)H
  \ar[dd]|-{(\mathrm{sw})(\id)(\id)}
  \\
  ({}_0W^*)({}_{\phantom{*}0}^{**}M)({}_0W)({\mathcal{I}}(N^{**})_0)H 
  \ar@{=}[ur]\ar[dd]_{\id\hat\tau_{N^{**}}}
  &
  ({\mathcal{I}}(N^{**})_0)({}_0W^*)({}_{\phantom{*}0}^{**}M)({}_0W) H
  \ar[d]_{\id\hat\nu_M^{-1}}
  &
  \\
  &
  ({\mathcal{I}}(N^{**})_0)({}_0W^*)({\mathcal{I}}(M^{**})_0)({}_0W) H
  \ar[d]\ar[r]
  &
  ({\mathcal{I}}((MN)^{**})_0)({}_0W^*)({}_0W)H
  \ar[dd]^{(\id)\mathrm{ev}(\id)}
  \\
  ({}_0W)({}^{**}_{\phantom{*}0}M)({}_0W)({}_0N^{**})H
  \ar[dddd]_{({}_0\mu_M)(\id)(\id)}
  &
  ({\mathcal{I}}(N^{**})_0)({\mathcal{I}}(M^{**})_0)({}_0W^*)({}_0W) H
  \ar[d]_{(\id)(\id)\mathrm{ev}(\id)}
  &
  \\
  &
  ({\mathcal{I}}(N^{**})_0)({\mathcal{I}}(M^{**})_0)H
  \ar[d]_{(\id)\tau_{M^{**}}}\ar[r]^-\cong
  &
  ({\mathcal{I}}((MN)^{**})_0)H
  \ar[ddd]^{\tau_{(MN)^{**}}}
  \\
  &
  ({\mathcal{I}}(N^{**})_0)({}_0M^{**})H
  \ar[d]_{\mathrm{sw}(\id)}
  &
  \\
  &
  ({}_0M^{**})({\mathcal{I}}(N^{**})_0)H\ar[dl]^{\id\tau_{N^{**}}}
  &
  \\
  ({}_0M^{**})({}_0N^{**})H
  \ar[rr]_-\cong\ar@{}[uuuuur]|-{\scriptstyle\mathrm{(A)}}
  &&
  ({}_0(MN)^{**})H
  \ar@{}[uuul]|-{\scriptstyle\mathrm{(B)}}
  \enddiagram}
$$
  \caption{{}}
  \label{fig:mumon1}
\end{figure}
\begin{figure}
$$
\diagramcompileto{mumon2}
({\mathcal{I}}(M)_0\otimes H)\square_H({\mathcal{I}}(N)_0\otimes H)
\ar[r]^-\cong\ar[ddd]_{\tau_M\square_H\tau_N}\ar@{}[dddr]|-{\mathrm{(D)}}
&
{\mathcal{I}}(N)_0\otimes{\mathcal{I}}(M)_0\otimes H
\ar[r]\ar[d]^{\id\otimes\tau_M} \ar@{}[dddr]|-{\mathrm{(E)}}
&
{\mathcal{I}}(M\otimes N)_0\otimes H\ar[ddd]^{\tau_{M\otimes N}}\\
&{\mathcal{I}}(N)_0\otimes{}_0M\otimes H\ar[d]^{\mathrm{sw}\otimes\id}
&\\
&{}_0M\otimes{\mathcal{I}}(N)_0\otimes H\ar[d]^{\id\otimes\tau_N}&\\
({}_0M\otimes H)\square_H({}_0N\otimes H)\ar[r]^-\cong&
{}_0M\otimes{}_0N\otimes H\ar[r]&
{}_0(M\otimes N)\otimes H
\enddiagram
$$
\caption{}
\label{fig:mumon1.1}
\end{figure}
\begin{figure}
$$
\diagramcompileto{mumon3}
{}_0W^*\otimes{}^{**}_{\phantom{*}0}M\otimes_0W\otimes
{\mathcal{I}}(N^{**})_0\otimes H
\ar[rr]^-{\mathrm{sw}\otimes\id}
\ar[d]|-{\id\otimes\id\otimes\id\otimes\tau_{N^{**}}}
\ar[ddrr]|-{{}_0\mu_M\otimes\id\otimes\id}
&&
{\mathcal{I}}(N^{**})_0 \otimes
{}_0W^*\otimes{}^{**}_{\phantom{*}0}M\otimes_0W 
\otimes H\ar[d]^{\id\otimes {}_0\mu_M\otimes\id}\\
{}_0W^*\otimes{}^{**}_{\phantom{*}0}M\otimes{}_0W\otimes{}_0N^{**}\otimes
H
\ar[d]_{{}_0\mu_M\otimes\id\otimes \id}
&&
{\mathcal{I}}(N^{**})_0\otimes{}_0M^{**}\otimes H
\ar[d]^{\mathrm{sw}\otimes\id}
\\
{}_0M^{**}\otimes{}_0N^{**}\otimes H
&&
{}_0M^{**}\otimes{\mathcal{I}}(N^{**})_0\otimes H
\ar[ll]^{\id\otimes\tau_{N^{**}}}
\enddiagram
$$
\caption{}
\label{fig:mumon3}
\end{figure}
\begin{figure}
$$
\diagramcompileto{mumon3}
{}_{\phantom{*}0}^{**}M\otimes
{}_{\phantom{*}0}^{**}N\otimes{}_0W\otimes H
\ar[rr]^-{\cong}\ar[d]_{\id\otimes\nu^{-1}_N}
&&
{}_{\phantom{*}0}^{**}(M\otimes N)\otimes{}_0(M\otimes
N)\otimes{}_0W\otimes H
\ar[d]^{\nu_{M\otimes N}^{-1}}
\\
{}_{\phantom{*}0}^{**}M\otimes
{}^{**}{\mathcal{I}}(N)_0\otimes{}_0W\otimes H
\ar[d]_{\mathrm{sw}\otimes \id}
&&
{}^{**}{\mathcal{I}}((M\otimes N)^{**})_0\otimes{}_0W\otimes H
\\
{}^{**}{\mathcal{I}}(N)_0\otimes {}_{\phantom{*}0}^{**}M
\otimes{}_0W\otimes H
\ar[rr]_-{\id\otimes\nu_M^{-1}}
&&
{}^{**}{\mathcal{I}}(N)_0\otimes{}^{**}{\mathcal{I}}(M)_0\otimes{}_0W
\otimes H\ar[u]_-\cong
\enddiagram
$$  
  \caption{ }
  \label{fig:mumon4}
\end{figure}
\begin{figure}
$$
\xymatrixcolsep{1.8cm}
\diagramcompileto{mumon5}
{}^{**}_{\phantom{*}0}M\otimes{}^{**}_{\phantom{*}0}N{}^{**}H
\ar[r]^-\cong\ar[d]_{\id\otimes{}^{**}\tau_N^{-1}}
&
{}^**_{\phantom{*}0}(M\otimes N)\otimes{}^{**}H
\ar[d]^{{}^{**}\tau_{M\otimes N}^{-1}}
\\
{}^{**}_{\phantom{*}0}M\otimes{}^{**}{\mathcal{I}}(N)_0\otimes
{}^{**}H
\ar[d]_{\mathrm{sw}\otimes \id}
&
{}^{**}{\mathcal{I}}(M\otimes N)_0\otimes{}^{**}H
\\
{}^{**}{\mathcal{I}}(N)_0\otimes{}^{**}_{\phantom{*}_0}M\otimes{}^**H 
\ar[r]_-{\id\otimes{}^{**}\tau_M^{-1}}
&
{}^{**}{\mathcal{I}}(N)_0\otimes
{}^{**}{\mathcal{I}}(M)_0\otimes{}^{**}H 
\ar[u]_\cong
\enddiagram
$$  
\caption{ }
  \label{fig:mumon5}
\end{figure}

{\em First axiom.}
The image of the diagram on the left hand side of Figure
\ref{fig:mumon0} is the exterior rectangle in Figure \ref{fig:mumon1}.
So it is enough to show the latter commutes, as the functor $M\mapsto
{}_0M\otimes H$ is an equivalence by Corollary \ref{c:ft}. The sub
diagrams left blank commute trivially. 

The diagram marked by (A) is just the commutative rectangle in Figure
\ref{fig:mumon3}. This is easy to show using the naturality of
$\mathrm{sw}$ and the definition of $\mu$. 
The diagram (B) in Figure \ref{fig:mumon1} commutes if the diagram
marked by (E) in Figure \ref{fig:mumon1.1} commutes for all $M,N$. To
show this, observe that the exterior rectangle in Figure
\ref{fig:mumon1.1} commutes by monoidality of $\tau$ and that the sub
diagram (D) commutes by Observation \ref{o:fsquareg}.

Finally, the diagram marked by (C) in Figure \ref{fig:mumon1.1}
commutes if and only if the diagram in Figure \ref{fig:mumon4}
does. If we tensor this diagram with ${}_0W^*$ on the right, after
composing with the isomorphism $\gamma:{}_0W\otimes
H\otimes{}_0W^*\cong {}^{**}H$ of Lemma \ref{lema:gamma}, we get the
diagram in Figure \ref{fig:mumon5}, which commutes as the diagram (E)
referred to above does. 

{\em Second axiom.}
Now we prove the commutativity of the diagram involving $\Bbbk$ in Figure
\ref{fig:mumon0}. For this we will need some notation. The symbol
$\Bbbk$ will denote the trivial {\em left}\/ $H$-comodule.
Let us denote
the canonical isomorphisms between $\Bbbk$ and both ${}^{**} \Bbbk$ and
$\Bbbk^{**}$ by $j$. As ${\mathcal{I}}$ is a monoidal functor
--see Corollary \ref{coro:iotamonoidal}--, we have canonical isomorphisms
$\delta:{}^{**}{\mathcal{I}}(\Bbbk)\to {\mathcal{I}}(\Bbbk^{**})$ in
$\M^H$ and 
$\theta:{\mathcal{I}}(\Bbbk)_0\to {}_0\Bbbk$ in ${}^H\!\M^H$. One
useful observation is that, as the monoidal structure on
${\mathcal{I}}$ is the unique one making $\tau$ a monoidal natural
transformation, we have $\theta\otimes
\id=\tau_\Bbbk:{\mathcal{I}}(\Bbbk)_0\otimes H\to{}_0\Bbbk\otimes H$.

The diagram we want to show that commutes lies in the category
$\M^H$; hence, we may equivalently show that its image under the
functor $M\mapsto {}_0M\otimes H$ of Corollary \ref{c:ft}
commutes. This new diagram is the one outer diagram in Figure 
\ref{fig:mumon6}. Indeed, the image of $\mu_\Bbbk$ is the
arrow $\zeta$ in \eqref{eq:zeta2}, that, by naturality of
$\mathrm{sw}$, is equal to the composition
$(\delta\otimes\id)(\mathrm{ev}\otimes\id)
(\id\otimes\mathrm{sw}\otimes\id)(\id\otimes\nu^{-1}_{\Bbbk})$ on the
right hand side of the diagram in Figure \ref{fig:mumon6}. The only
sub diagrams whose commutativity is not obvious are the ones marked
with (a), (b) and (c). 

The commutativity of the diagram (a) follows form the following
observation. By definition, $\nu_\Bbbk^{-1}\otimes \id_{{}_0^*W}$
corresponds, up to composing with certain canonical isomorphisms,  
to ${}^{**}\tau_{\Bbbk}$ --see Lemma \ref{lema:xi}--. On the other
hand,
${}^{**}\theta\otimes\id\otimes\id\otimes\id:
{}^{**}{\mathcal{I}}(\Bbbk)_0\otimes{}_0W\otimes H\otimes {}_0W^*\to
{}_{\phantom{*}0}^{**}\Bbbk\otimes {}_0W\otimes H\otimes {}_0W^*$
  also corresponds, up to composing with the same canonical isomorphisms,
  to ${}^{**}\tau_\Bbbk$ (this because
  $\theta\otimes\id_H=\tau_\Bbbk$). It follows that the diagram (a)
  commutes. 

The diagram made by (b) commutes because $\theta$ is induced by the
monoidal structure of ${\mathcal{I}}$, and monoidal functors preserve
duals. 

Finally, the diagram (c) commutes by naturality of $\tau$. 

\begin{figure}
$$
\xymatrixcolsep{1.6cm}
\diagramcompileto{mumon6}
&({}_0W^*)({}^{**}_{\phantom{*}0}\Bbbk)({}_0W)H
\ar[dr]^{\id\nu^{-1}_\Bbbk}
\ar@{}@<25pt>[d]|-{\mathrm{(a)}}
&
\\
({}_0W^*)({}_0\Bbbk)({}_0W)H
\ar[r]_-{(\id) j(\id)}\ar[ur]^{(\id) j(\id)(\id)}
&
({}_0W^*)({}_{\phantom{*}0}^{**}\Bbbk)({}_0W)
\ar@{=}[u]\ar[d]|-{(\id)\mathrm{sw}(\id)}
&
({}_0W^*)({}{**}{\mathcal{I}}(\Bbbk)_0)({}_0W)H
\ar[d]^{(\id)\mathrm{sw}(\id)}
\ar[l]^-{(\id){}^{**}\theta(\id)(\id)}
\\
&({}_0W^*)({}_0W)({}^{**}_{\phantom{*}0}\Bbbk)H
\ar[d]_{{\mathrm{ev}}(\id)}
&
({}_0W^*)({}_0W)({}^{**}{\mathcal{I}}(\Bbbk)_0)H
\ar[d]^{(\mathrm{ev})(\id)}
\ar[l]^-{(\id)(\id){}^{**}\theta(\id)}
\\
({}_0\Bbbk)H
\ar[uu]^\cong\ar[r]^-{j\id}
\ar[dr]|-{\tau_\Bbbk=\theta(\id)}\ar[ddr]_{j\id}
&
({}_{\phantom{*}0}^{**}\Bbbk)H
\ar[dr]_\cong\ar@{}[d]|-{\mathrm{(b)}}
&
({}^{**}{\mathcal{I}}(\Bbbk)_0)H
\ar[l]_-{{}^{**}\theta(\id)}
\ar[d]^{\delta(\id)}
\\
&({\mathcal{I}}(\Bbbk)_0)H
\ar[r]_-{({\mathcal{I}}(j)_0)(\id)}
\ar@{}[d]|-{{\mathrm{(c)}}}
&
({\mathcal{I}}(\Bbbk^{**})_0)H\ar[dl]^{\tau_{\Bbbk^{**}}}
\\
&({}_0\Bbbk^{**})H
\enddiagram
$$  
  \caption{ }
  \label{fig:mumon6}
\end{figure}
\end{proof}

The monoidality of the natural transformation $\mu$ just proved
translates into properties of the functional $\sigma$ in Theorem
\ref{theo:radford}. We compute the below in an explicit way the
monoidal structure of $\sigma$.

\begin{obse} The functional $\sigma$ induces the
natural isomorphism of Corollary \ref{c:Msquare()+}, which is monoidal
since $\mu$ is. Therefore, if we know the monoidal structures of the
morphisms $p(a^{-1}\otimes p(\bar S^2\otimes a))$ and $S^2$, we can
deduce the equations satisfied by $\sigma$. More explicitly, if
these morphisms have monoidal structures $(\chi_1,\rho_1)$ and
$(\chi_2,\rho_2)$ respectively, then $\sigma$ satisfies 
\begin{equation}
  \label{eq:sigmapropr1}
  \chi_1\star \sigma p=(\sigma\otimes\sigma)\star \chi_2\qquad
  \rho_1\sigma(1)=\rho_2.
\end{equation}

The antipode $S:H^\circ\to H$ has a monoidal structure $(\chi^S,1)$,
where $\chi^S$ is given explicitly in Proposition 
\ref{obs:dualsrlcirc2}. By Observation \ref{o:compmonmorph} we have
that $S^2$, which is the composition of $S^{\mathrm{cop}}:H\to
H^{\circ}$ and $S:H^\circ\to H$, has $(\chi^S(S\otimes S)\star
(\chi^S)^{-1}\mathrm{sw},1)$ as monoidal structure. This is because
$S^{\mathrm{cop}}$ has a monoidal structure $((\chi^S)^{-1}\mathrm{sw},1)$. 
The inverse of the antipode $\bar S:H^\circ\to H$ has a canonical
monoidal structure given in terms of $\chi^S$ by $((\chi^S)^{-1}(\bar
S\otimes\bar S),1)$. Thus, $\bar S^2$, this is, the composition of
$\bar S$ with $\bar S^{\mathrm{cop}}:H^\circ\to H$, has a monoidal
structure $(\chi^S(\bar S^2\otimes \bar S^2)\mathrm{sw}\star
(\chi^S)^{-1}(\bar S\otimes \bar S),1)$.  

The morphism $p(a^{-1}\otimes p(\id\otimes a)):H\to H$ is monoidal
with a monoidal structure given by $(\chi_0,1)$ where $\chi_0$ is the
following product in $(H\otimes H)^\vee$. 
$$
\phi^{-1}(a^{-1}\otimes (-) a \otimes a^{-1}((?)a))\star
\phi((-)a\otimes a^{-1}\otimes (?)a)\star
\phi^{-1}(-\otimes a\otimes a^{-1})\star
\phi(-\otimes ?\otimes a)
$$
Then, the monoidal structure $(\chi_1,\rho_1)$ of the composition of
$\bar S^2$ with $p(a^{-1}\otimes p(\id\otimes a))$ is given by
$\chi_1=\chi_0(\bar S^2\otimes \bar S^2)\star \chi^S(\bar S^2\otimes
\bar S^2)\mathrm{sw}\star(\chi^S)^{-1}(\bar S\otimes \bar S)$ 
and $\rho_1=1$. We deduce that $\sigma$ satisfies $\sigma(1)=1$ and
\begin{equation}
  \label{eq:sigmapropr2}
  \chi_0(\bar S^2\otimes \bar S^2)\star \chi^S(\bar S^2\otimes
\bar S^2)\mathrm{sw}\star (\chi^S)^{-1}(\bar S\otimes\bar S)
\star \sigma p=
(\sigma\otimes \sigma)\star
\chi^S(S\otimes S)\star(\chi^S)^{-1}\mathrm{sw}.
\end{equation}
\end{obse}
\section{The case of a Hopf algebra}\label{s:casehopf}
We briefly mention the needed adjustments to the proof above in order
to get the classical Radford's formula for $S^4$. We assume that $H$
is a finite dimensional Hopf algebra and define the following functions.
 
Denote by  $\omega \in H^{\vee}$, the modular function of $H$ that we
know it 
is an algebra homomorphism.  It
can be defined as the modular element in the Hopf algebra $H^{\vee}$ 
(the linear dual of
$H$). In particular if $i \in H$ is a right integral, the functional
$\omega$ is characterized by the property that for all $x \in H$ we
have that $xi=\omega(x)i$. 

We will also consider the automorphism of Nakayama $\mathcal N$, that is
characterized by the equation $\varphi(xy)=\varphi(y \mathcal N(x))$ for all
$x,y \in H$ where $\varphi$ is as before a right cointegral for $H$. 

Next we show how to obtain an expression of 
the inverse of Nakayama's automorphism
in terms of $S$ and $\omega$.   
For all $x,y \in H$, we have that:
\begin{equation}\label{e:fourierhopf} 
\sum \varphi(y_1 x)y_2= \sum \varphi(y_1x_1) y_2\varepsilon(x_2) = \sum
\varphi(y_1x_1)y_2x_2 S(x_3)=\sum \varphi(y
  x_1)S(x_2)
\end{equation} 
If we take 
  $y=i$ in the above equality and assume that $\varphi(i)=1$ 
we obtain that $S(x)=\sum \varphi(i_1x)i_2$
  or equivalently that $x=\sum \varphi(i_1x)\overline{S}i_2$.
Hence, it follows that $\mathcal
N(x)=\sum \varphi(i_1\mathcal Nx)\overline{S}i_2= \overline{S}^2(\sum \varphi(xi_1)Si_2)$. 
Now, using using again the equation \eqref{e:fourierhopf} 
we conclude that 
\begin{equation*}
\mathcal N(x)=\overline{S}^2(\sum
\varphi(x_1i)x_2)=\overline{S}^2(\sum \omega(x_1)x_2)=
\overline{S}^2(x
\leftharpoonup \omega).
\end{equation*} 
Taking the inverse maps in the above equation
we conclude that $S^2x \leftharpoonup \omega^{-1} =\mathcal
N^{-1}x$. Then, it follows that  $\varepsilon
\mathcal N^{-1} =\omega^{-1}$.

In the case of a Hopf algebra, for $M \in \rcomodH$ then 
  $\tau_M: M \otimes H \rightarrow M \otimes H$ is given by the
  formula $\tau_M(m \otimes h) = \sum m_0 \otimes S(m_{1})h$. This
  is easily obtained from Theorem \ref{lema:tau} substituting the
  associators as well as $\alpha$ and $\beta$ by $\varepsilon$.   
  The inverse of $\tau_M$ is given as $\tau_M^{-1}(m \otimes h) = 
\sum m_0 \otimes  m_{1}h$.  

With respect to the properties of duality, the same formul\ae\/
\eqref{e:M*structure} and \eqref{e:*Mstructure} yields the comodule
structure on the dual spaces. The evaluation and coevaluation in this
case are the same than the usual ones in the category of vector
spaces (see formul\ae \/ \eqref{eq:ev}, \eqref{eq:coev}, \eqref{eq:evr}, 
\eqref{eq:coevr}). 

If $M$ is a left $H$--module the natural 
right $H$--module structure on the left dual --see \eqref{e:M*structure}--
$a_{{}^*M}: {}^{*}M \otimes H \rightarrow {}^*M$ is given by
$$
{}^*M\otimes H\xrightarrow{\id \otimes \id \otimes\mathrm{c}}{}^*M\otimes
H\otimes M\otimes {}^*M\xrightarrow{\id\otimes a_M\otimes
  \id}{}^*M\otimes M\otimes {}^*M\xrightarrow{\mathrm{e}\otimes \id}{}^*M.
$$
For the right dual the formula is similar --see \eqref{e:M*structure}.
In particular, in the case we consider $H \in \bicomodH$ as a left
module with respect to the regular action, the right 
$H$--structure considered above in
this situation is simply the following: $ f \leftharpoonup h \in
{}^*H$, $(f \leftharpoonup h)(x) = f(hx)$. 

The Frobenius map $\mathcal F$, that was given in Theorem \ref{lema:basic}, is 
$\mathcal F(\varphi \otimes h)= \varphi
\leftharpoonup h$ for $\varphi \in W$ and $h \in H$. 

In particular as we mentioned in Observation
\ref{obse:generalgamma} part b), 
the morphism of algebras $\gamma$ defined in Lemma \ref{lema:gamma}:
$$\gamma: {}_0 W \otimes H \otimes {}_0{}^* W
  \xrightarrow{\mathcal F \otimes \id}{}^*H \otimes {}_0{}^* W \cong
  {}^* ({}_0 W \otimes H)\xrightarrow{({}^*\mathcal{F})^{-1}}
  {}^{**}H\/,$$ is given by the formula: $\gamma(\varphi \otimes x
  \otimes t)= \mathrm{ev}(- \otimes {\mathcal N}x)$ where
$\mathcal N: H \rightarrow H$ is the Nakayama morphism, that in this
case is an algebra automorphism.      

The map
    $\xi_M$ considered in Lemma \ref{lema:xi}, can be described
    explicitly by 
$\xi_M(\mathrm{ev}_m \otimes \varphi \otimes h \otimes t)= \sum
    \mathrm{ev}_{m_0}\otimes \varphi \otimes\mathcal N^{-1}(S(m_{1}))h
    \otimes t$.
 Indeed, from the commutative diagram \eqref{e:xi} 
we deduce that $\omega_M(\mathrm{ev}_m \otimes
    \mathrm{ev}_h)= \sum \mathrm{ev}_{m_0}\otimes
    \mathrm{ev}_{S(m_{1})h}$. To prove the formula for $\xi_M$ 
we prove that 
\begin{equation*}
(\id \otimes
    \gamma)(\sum
    \mathrm{ev}_{m_0}\otimes \varphi \otimes\mathcal N^{-1}(S(m_{1}))h
    \otimes t)= \omega_M(\id \otimes \gamma)(\mathrm{ev}_m
    \otimes \varphi \otimes h \otimes t).
\end{equation*}
The left hand side of the above equation is: 
\begin{equation*}
(\id \otimes \gamma) (\sum \mathrm{ev}_{m_0}\otimes \varphi \otimes
\mathcal N^{-1}(S(m_{1}))h \otimes t)
= \sum \mathrm{ev}_{m_0} \otimes \mathrm{ev}_{S(m_{1})\mathcal N(h)},
\end{equation*}
while the right hand side can be computed as:
\begin{equation*}
\omega_M(\id \otimes \gamma)(\mathrm{ev}_m \otimes \varphi
    \otimes h \otimes t)
=\omega_M(\mathrm{ev}_{m} \otimes \mathrm{ev}_{\mathcal N(h)})
= \sum\mathrm{ev}_{m_0}\otimes \mathrm{ev}_{S(m_{1})\mathcal N(h)}.
\end{equation*}

Hence the map $\nu_M: {}^{**}{\mathcal{I}}(M)_0\otimes {}_0W\otimes H \to{}^{**}
  _{\hspace{.1cm}0}M\otimes {}_0W\otimes H$ introduced in Definition
  \ref{defi:nueme} is given by: $\nu_M(\mathrm{ev}_m \otimes \varphi \otimes h)= \sum
    \mathrm{ev}_{m_0}\otimes \varphi \otimes\mathcal N^{-1}(S(m_{1}))h$.
Moreover, the map $\widehat \nu_M$ has exactly the same expression than
$\nu_M$. 

For later use we record the following formula for  $\widehat
\nu_M^{-1}$ that can be proved by a direct computation:
\begin{equation}\label{eq:nuhatinv}
\widehat \nu_M^{-1}(\operatorname {ev}_m \otimes \varphi
\otimes h)= \sum \operatorname{ev}_{m_0} \otimes \varphi \otimes
\mathcal N^{-1}(m_1)h.
\end{equation}

Thus, the morphism $\zeta_M$
defined in \eqref{eq:zeta2} is given as: 
$$\zeta_M(t \otimes \mathrm{ev}_m \otimes \varphi \otimes h)= \sum
\mathrm{ev}_{m_0} \otimes 
\overline{S}(m_1) \mathcal N^{-1}(m_{2})h.
$$
Indeed it follows from equation \eqref{e:M*structure} that the right coaction
  in $M^{**}$ is $\chi_{M^{**}}(\mathrm{ev}_m)= \sum 
\mathrm{ev}_{m_0} \otimes \overline{S}^2(m_1)$. Thus, 
\begin{multline*}
\zeta_M(t \otimes \mathrm{ev}_m \otimes \varphi \otimes h)
=(\mathrm{ev} \otimes \tau_{M^{**}})(\id \otimes \mathrm{sw}\otimes
\id)(\id \otimes \widehat \nu_M^{-1})
(t \otimes \mathrm{ev}_m \otimes \varphi \otimes h)= \\ 
\sum (\mathrm{ev}
\otimes \tau_{M^{**}})(\id \otimes \mathrm{sw}\otimes \id)(t \otimes
\mathrm{ev}_{m_0} \otimes \varphi \otimes \mathcal N^{-1}(m_1)h)
= \\ 
\tau_{M^{**}}(\sum \mathrm{ev}_{m_{0}}  \otimes \mathcal N^{-1}(m_1)h)
=  \sum \mathrm{ev}_{m_{0}}  \otimes
\overline{S}(m_1)\mathcal
N^{-1}(m_{2})h. 
\end{multline*}

Then, as $\zeta_M=\mu_M \otimes
  \id_H$, with $\mu_M: W^* \otimes{} ^{**}M \otimes W \rightarrow
  M^{**}$, it is clear that $\mu_M$ satisfies the following
  equality: 
$\mu_M(t \otimes \mathrm{ev}_m \otimes \varphi) \otimes h = \sum
  \mathrm{ev}_{m_{0}}  \otimes \overline{S}(m_1)\mathcal N^{-1}(m_{2})h$.
 
If we apply $\id  \otimes \id  \otimes \varepsilon$
  to the equality above we obtain:  
$$\mu_M(t \otimes \mathrm{ev}_m \otimes \varphi) = \sum
  \mathrm{ev}_{m_{0}} (\varepsilon\mathcal
  N^{-1})(m_{1})= \sum
  \mathrm{ev}_{m_{0}} \omega^{-1}(m_{1})$$
We have used above the equality 
$\varepsilon \mathcal N^{-1}= \omega^{-1}$ proved
before. Hence we deduce that 
$\mu_M(t \otimes \mathrm{ev}_m \otimes \varphi) = 
  \mathrm{ev}_{\omega^{-1} \rightharpoonup m}$.

Next, we observe that the natural isomorphism  constructed
in Corollary \ref{c:Msquare()+} is simply the map $m \mapsto
(\omega^{-1}\rightharpoonup m)$. 
Applying the bijections proved in Theorem
\ref{l:nacho}, we find that the map $\sigma$ appearing in Radford's
formula --Theorem \ref{theo:radford}-- is simply
$\sigma(h)=\varepsilon (\omega^{-1} \rightharpoonup h)= \omega^{-1}(h)$.   
Hence, we deduce the classical Radford's formula
$$a^{-1}\overline{S}^2(x)a = \omega^{-1} \rightharpoonup S^2(x)
\leftharpoonup \omega \quad\text{or}\quad S^4(x) = \omega
\rightharpoonup a^{-1} x a \leftharpoonup \omega^{-1}.
$$ 
This shows that the functional $\sigma$ is indeed the coquasi Hopf
algebra analogue of the modular function. 

Next we explain how the monoidality of $\sigma$ proved at the end of
the previous section generalizes the multiplicativity of the modular
function $\omega\in H^\vee$. 

Recall that in the case of a Hopf algebra the associativity of the
product and the fact that $S$ is a morphism of algebras are expressed
as $\phi=\varepsilon\otimes\varepsilon\otimes\varepsilon$ and
$\chi^S=\varepsilon\otimes\varepsilon$. Therefore, the equality
\eqref{eq:sigmapropr2} simplifies to $\sigma p=\sigma\otimes\sigma$,
this is, $\sigma$ is multiplicative. The equality $\sigma(1)=1$ was
shown for an arbitrary coquasi Hopf algebra. Hence
$\omega=\sigma^{-1}$ is a morphism of algebras. 

An important point is that in the proof
of the monoidality of $\sigma$ we did not use integrals (only
cointegrals, presented as the comodule $W$). Then, in the case of a
finite dimensional Hopf algebra, if we use $\sigma$ instead of the
modular function, we can avoid mentioning integrals altogether in the
proof of the classical Radford's formula. 

\section{Appendix: categorical background}\label{s:back}
This appendix is an account of some basic results on
functors between categories of comodules. These results are not
unknown, but the proofs found in the literature are usually {\em ad
  hoc}.  The unified presentation below is based on density of functors and
completions of categories under certain classes of colimits. 

We will work with categories enriched in the category of vector
spaces over a field $\Bbbk$, sometimes called $\Bbbk$--linear
categories. 
Although one has to be careful when dealing with enriched 
categories, in our case the subtleties of the theory disappear. This
is a consequence of the fact that the {\em underlying set functor}\/
$\mathbf{Vect}(\Bbbk,-):\mathbf{Vect}\to\mathbf{Set}$ is conservative
({\em i.e.}, reflects isomorphisms). We denote by $[\mathcal
A,\mathcal B]$ the $\Bbbk$--linear category of $\Bbbk$--linear
functors $\mathcal A\to\mathcal B$ and natural transformations between
them.

Recall the notion of dense functor (see \cite{kn:kelly} for a complete
exposition on the subject). Let $\mathcal A$ be a small category. 
A functor $K:\mathcal A\to \mathcal C$ is
dense if the functor $\tilde K:\mathcal C\to [\mathcal
A^{\mathrm{op}},\mathbf{Vect}]$, given by $C\mapsto \mathcal C(K-,C)$,
is fully faithful. In other words, $K$ is dense if every natural
transformation $\mathcal C(K-,C)\Rightarrow \mathcal C(K-,D)$ is of
the form $\mathcal C(K-,f)$ for a unique $f:C\to D$ in $\mathcal C$. 
When $K$ is the inclusion of a full
subcategory, say that $\mathcal A$ is dense in $\mathcal C$. 

A colimit in $\mathcal C$ is {\em $K$-absolute}\/ if it is preserved
by $\tilde K$. In elementary terms, a colimit $\sigma_j:P(j)\to
\operatorname{colim} P$ of a functor $P:\mathcal J\to \mathcal C$ is
$K$-absolute if for all objects $A\in \mathcal A$ the transformation
$\mathcal C(K(A),P(j))\to\mathcal C(K(A),\operatorname{colim}P)$ is a
colimit in $\mathbf{Vect}$. 

Now suppose $K$ is the inclusion of a full subcategory $\mathcal A$
into $\mathcal C$. Consider a family of functors
$\Phi=\{P_\gamma:\mathcal J_\gamma\to\mathcal
C\}_{\gamma\in\Gamma}$. We say that $\mathcal C$ is the closure of
$\mathcal A$ under the family $\Phi$ if there is no proper
replete full subcategory $\mathcal D$ of $\mathcal C$ containing (the
image of)
$\mathcal A$ such that
$\operatorname{colim}P_\gamma\in\mathcal D$ whenever $P_\gamma$ takes
factors through $\mathcal D$. If $\Phi$ is the family of all
the functors with small (finite) domine into $\mathcal C$, we say that
$\mathcal C$ is the closure of $\mathcal A$ under small (finite)
colimits. 

We say that $\Phi$ is a {\em density presentation}\/ for
$K$ if the colimit of each
$P_\gamma$ exists and is $K$-absolute, and $\mathcal C$ is the closure
of $\mathcal A$ under the family $\Phi$. 
The functor $K$ is {\em dense}\/ when it has a density presentation
(although density can be defined in other ways, our choice here is
justified by \cite[Theorem 5.35]{kn:kelly}).

Write $\mathrm{Cocts}^K[\mathcal C,\mathcal B]$ for the full
sub-$\Bbbk$-category of $[\mathcal C,\mathcal B]$ of those functors
that preserve $K$-absolute colimits. 
The following is a particular instance of \cite[Thm. 5.31]{kn:kelly}.
\begin{theo}\label{t:dense}
Let   $\Phi=\{P_\gamma:\mathcal J_\gamma\to\mathcal
C\}_{\gamma\in\Gamma}$ be a density prsentation of the fully faithful
functor $K:\mathcal A\to\mathcal C$. Suppose each $\mathcal J_\gamma$
is small and that $\mathcal B$ admits all small colimits. Then
precomposing with $K$ yields an equivalence 
$$
\mathrm{Cocts}^K[\mathcal C,\mathcal B]\simeq[\mathcal A,\mathcal B]
$$
with pseudoinverse given by taking left Kan extensions along $K$. 
\end{theo}

A basic example of dense subcategory is provided by the category of
modules over a ring $R$. If we take $\mathcal C={}_R\M$
and $\mathcal A$ the full subcategory
determined by the $R$-module $R$, then $\mathcal A$ is dense in
$\mathcal C$. A density presentation is given by the family of
functors with small domain into ${}_R\M$. This is so
because all colimits are $K$-absolute: after identifying $[\mathcal
A^{\mathrm{op}},\mathbf{Vect}]$ with ${}_R\M$, $K$ is isomorphic to
the identity functor. 
The category $\mathcal A$ is also dense in the category
${}_R\M_{\mathrm{fp}}$ of finitely presented
$R$-modules. Indeed, ${}_R\M_{\mathrm{fp}}$ is the
closure of $\mathcal A$ under finite colimits and these are
$K$-absolute, where $K$ is the inclusion of $\mathcal A$ into
${}_R\M_{\mathrm{fp}}$. Observe that $\tilde K$ in this
case is isomorphic to the inclusion of
${}_R\M_{\mathrm{fp}}$ into ${}_R\M$, and
hence it preserves colimits. As a consequence of Theorem \ref{t:dense}
we have equivalences
$$
\mathrm{Cocts}[{}_R\M,\mathcal B]\simeq [\mathcal
A,\mathcal B] = R^{\mathrm{op}}\text-\mathcal B \qquad 
\mathrm{Rex}[{}_R\M_{\mathrm{fp}},\mathcal
D]\simeq[\mathcal A,\mathcal D]=R^{\mathrm{op}}\text-\mathcal D
$$
for any categories $\mathcal B$ and $\mathcal D$ with small colimits
and finite colimits respectively. Here $R^{\mathrm{op}}\text-\mathcal
B$ denotes the 
category with objects $B$ of $\mathcal B$ equipped with an action or
$R^{\mathrm{op}}$, that is, a ring morphism
$R^{\mathrm{op}}\to\mathcal B(B,B)$, and evident morphisms. 

Slightly more general, the Yoneda embedding $K:\mathcal A\to[\mathcal
A^{\mathrm{op}},\mathbf{Vect}]$ is dense for any small category
$\mathcal A$. 

A second example of interest for us is subcategory of
finite-dimensional comodules $\M^C_f$. Let $K:\M^C_f\to\M^C$ be the
inclusion functor. Given a
$C$-comodule $M$, consider the {\em comma category}\/ $K/M$. That is,
the category whose objects are pairs $(N,f)$ where $N\in\M^C_f$ and
$f:N\to M$, and whose arrows $(N,f)\to(N',f')$ are the arrows $g:N\to
N'$ such that $f'g=f$. The functor $P_M:K/M\to \M^C$ sending $(N,f)$
to its $N$ is the base of a cone of
vertex $M$,  with components $\sigma_{(N,f)}=f:N \to
M$. Clearly $K/M$ is small and filtered (since $\M^C_f$ has finite colimits and
$K$ preserves them) and  $\sigma$ is a colimiting cone. The family of
functors $P_M$ with $M$ a $C$-comodule is a density presentation for
$K$: clearly $\M^C$ is the closure of $\M^C_f$ under filtered
colimits and filtered colimits are preserved by $\tilde K$ since
finite dimensional comodules are finitely presentable ({\em i.e.},
$\M^C(N,-)$ preserves filtered colimits whenever $N$ is
finite-dimensional). Since $K$ preserves finite colimits, it is clear
that the image of $\tilde
K:\M^C\to[(\M^C_f)^{\mathrm{op}},\mathbf{Vect}]$ lies in the full
subcategory $\mathrm{Lex}[(\M^C_f)^{\mathrm{op}},\mathbf{Vect}]$ of
left exact functors; moreover, the replete image of $\tilde K$  can be
shown to be exactly this subcategory. This yields an equivalence 
$$
\mathrm{Fin}[\M^C,\mathcal B]\simeq [\M^C_f,\mathcal B]
$$
for any category $\mathcal B$ with filtered colimits, where the
category on the left hand side is the category of finitary ({\em
  i.e.}, filtered colimit-preserving) functors. 

Our next example is the 
full subcategory $\mathcal A$ of $\M^C$
determined by the regular comodule $C$. Recall that each comodule $M$
is the equalizer of the canonical pair of comodule morphisms between
free comodules $M\otimes
C\to M\otimes C\otimes C$.  Clearly, $\M^C$ is the closure
of $\mathcal A$ under small limits. For, every free comodule has to be
in the closure under small limits, and hence each comodule does
too. Now consider the inclusion functor $K:\mathcal
A^{\mathrm{op}}\hookrightarrow (\M^C)^{\mathrm{op}}$. 
We shall show that the functor $\tilde K:(\M^C)^{\mathrm{op}}\to[\mathcal
A,\mathbf{Vect}]$ preserves small colimits. To do this, it is enough
to show that the composition of $\tilde K$ with the ``forgetful''
functor $[\mathcal A,\mathbf{Vect}]\to\mathbf{Vect}$ (recall that
$\mathcal A$ has just one object) preserves small colimits. The
resulting functor $(\M^C)^{\mathrm{op}}\to\mathbf{Vect}$ is simply
$\M^C(-,C)$, which is isomorphic to $\mathbf{Vect}(U(-),\Bbbk)$, where
$U$ denotes the forgetful functor $\M^C\to\mathbf{Vect}$. Clearly $U$
preserves limits and $\mathbf{Vect}(-,\Bbbk)$ converts them into
colimits. We have shown, then, that small colimits in
$(\M^C)^{\mathrm{op}}$ are $K$-absolute. This, together with the fact
that $(\M^C)^{\mathrm{op}}$ is the completion of $\mathcal
A^{\mathrm{op}}$ under small colimits, shows that $K$ is dense. We
get, then, an equivalence 
$$
\mathrm{Cts}[\M^C,\mathcal B]\simeq [\mathcal A,\mathcal B]
$$
for any category $\mathcal B$ with small limits; here the category on
the left hand side is the category of continuous ({\em i.e.}
small limit-preserving) functors and transformations between them. 

Finally, for a finite dimensional coalgebra $C$, 
consider the full subcategory $\mathcal A$ of $\M^C$ given by
the single object $C$. As $(\M^C_f)^{\mathrm{op}}$ is equivalent to
$(\M_{C^\vee})_f$, the functors $\mathcal J\to \mathcal
A^{\mathrm{op}}$ form a finite category $\mathcal J$ constitute a
density presentation for $\mathcal
A^{\mathrm{op}}\hookrightarrow(\M^C_f)^{\mathrm{op}}$. 
Furthermore, we have equivalences
$$
\mathrm{Lex}[\M^C_f,\mathcal
B]\cong\mathrm{Rex}[(\M^C_f)^\mathrm{op},\mathcal
B^{\mathrm{op}}]\simeq [\mathcal A^{\mathrm{op}},\mathcal
B^{\mathrm{op}}]\cong[\mathcal A,\mathcal B]
$$
for any category $\mathcal B$ with finite limits. 
Recall that there is an equivalence
$[\mathcal A,\mathbf{Vect}]\cong \M_{C^\vee}$. For, to give a
functor $\mathcal A\to \mathbf{Vect}$ is to give a vector space with
a left action of the algebra $\M^C(C,C)$, which is isomorphic to
$(C^\vee)^{\mathrm{op}}$ via the map sending $\gamma:C\to\Bbbk$ to
$(\gamma\otimes\mathrm{id})\Delta$. Now it is easy to deduce that for
a finite-dimensional coalgebra $C$ 
there are equivalences $\mathrm{Lex}[\M^C_f,\M^D]\simeq {}^C\!\M^D$,
sending a right exact functor $F$ to the bicomodule $F(C)$. A
pseudoinverse for this equivalence is the functor sending a bicomodule
$M$ to $-\square_C M$. 

In Section \ref{s:catsofbicomod} we used the following easy
observation. 

\begin{obse}\label{o:apendix}
Let $C$ be a finite dimensional coalgebra, $K:\M^C_f\to\M^C$ be the
inclusion functor and $M,N\in {}^C\!\M^D$. 
Then we have a string of canonical isomorphisms
\begin{align*}
{}^C\!\M^D(M,N)
&\cong\operatorname{Lex}[\M^C_f,\M^D](K(-)\square_CM,K(-)\square_CN)\\
&=[\M^C_f,\M^D](K(-)\square_CM,K(-)\square_CN)\\
&\cong \operatorname{Fin}[\M^C,\M^D](-\square_CM,-\square_CN)\\
&=[\M^C,\M^D](-\square_CM,-\square_CN)\\
\end{align*}
\end{obse}

The last piece of categorical background we will need is the {\em
  tensor product}\/ of categories with finite limits. This is closely
related to Deligne's tensor product of abelian categories of 
\cite{kn:deligne}. We only need the case of categories of finite
dimensional comodules over finite dimensional coalgebras, 
though, and in this case the existence of this
product reduces to few simple observations. 

Recall that the category $\mathcal C\otimes\mathcal D$ has as
objects pairs of objects $(c,d)$ with $c\in\mathcal C$ and
$d\in\mathcal D$, and homs $\mathcal C\otimes\mathcal
D((c,d),(c',d'))=\mathcal C(c,c')\otimes_\Bbbk\mathcal D(d,d')$. 
If $\mathcal C,\mathcal D,\mathcal E$ are categories with finite
limits, a functor $F:\mathcal C\otimes\mathcal D\to\mathcal E$ is {\em
  left exact in each variable}\/ if for each $c\in\mathcal C$ and
$d\in \mathcal D$ the functors $F(c,-):\mathcal D\to\mathcal E$ and
$F(-,d):\mathcal D\to\mathcal E$ are left exact. These functors,
together with the natural transformations between them, form a
category $\mathrm{Lex}[\mathcal C,\mathcal D;\mathcal E]$. 

A tensor
product of $\mathcal C$ with $\mathcal D$ as categories with finite
limits is a category with finite limits $\mathcal C\boxtimes\mathcal
D$ together with a functor $\mathcal C\otimes\mathcal D\to\mathcal
C\boxtimes\mathcal E$ left exact in each variable that induces
equivalences $\mathrm{Lex}[\mathcal C\boxtimes\mathcal D,\mathcal
E]\simeq\mathrm{Lex}[\mathcal C,\mathcal D;\mathcal E]$ for each
$\mathcal E$.  

The case of interest for us in this work is the one of categories of
finite dimensional comodules over finite dimensional coalgebras, dual
to the case considered in \cite{kn:deligne}. If $C,D$ are finite
dimensional coalgebras, we claim that the functor
$\otimes_\Bbbk:\M^C_f\otimes\M^D_f\to\M^{C\otimes D}_f$ is a tensor
product of $\M^C_f$ with $\M^D_f$ as categories with finite limits. To
see this, let us call $\mathcal C\subset\M^C_f,\mathcal
D\subset\M^D_f$ and $\mathcal B\subset\M^{C\otimes D}_f$ 
the full subcategories determined by the respective
regular comodule, and observe that there is a commutative diagram as
depicted below.
$$
\xymatrixrowsep{.6cm}
\diagramcompileto{productofcatscomod}
{\mathrm{Lex}}[{\M}^{C\otimes D}_f,{\mathcal{E}}]\ar[r]\ar[d]_\simeq&
{{\mathrm{Lex}}}[{\M}^C_f,{\M}^D_f;{\mathcal{E}}]\ar[d]^\simeq \\
[{\mathcal{B}},{\mathcal{E}}]\ar[r] & 
{[{\mathcal{C}}\otimes{\mathcal{{D}}},{\mathcal{E}}]}
\enddiagram
$$
The horizontal arrows are induced by $\otimes_\Bbbk$ and the obvious
functor $\mathcal C\otimes\mathcal D\to\mathcal B$ which at the level
of the unique hom-space is just
$\M^C_f(C,C)\otimes\M^D_f(D,D)\to\M^{C\otimes D}_f(C\otimes D,C\otimes
D)$. This last linear transformation is an isomorphism, as a
consequence of the finiteness of $C$ and $D$, and then the bottom row of
the diagram is an isomorphism. It follows that the top row is an
equivalence.


\begin{thebibliography}{99}

\bibitem{kn:Benabou} B{\'e}nabou, J. {\em  {Introduction to
      bicategories.}} {Reports of the Midwest Category Seminar.}
  Springer, Berlin. (1967) {pp. 1--77}.
\bibitem{kn:bkl} Bespalov, Y., Kerler, T., Lyubashenko, V. and Turaev,
  V. {\em Integrals for braided Hopf algebras.} J. Pure and Applied
  Algebra, {\bf 148}, (2000), pp. 113-164.  
\bibitem{kn:bulacu}
Bulacu D., Chiri{{t}}{\u{a}} B. {\em Dual {D}rinfeld double by diagonal
crossed coproduct}, Rev. Roumaine Math. Pures Appl. 47~(3) (2002)
271--294 (2003).
\bibitem{kn:bc} Bulacu, D. and Caenepeel, S. {\em Integrals for (dual)
    quasi--Hopf algebras, applications.}\,J. Algebra, {\bf 266}, 2, (2003),
pp. 552-583.
\bibitem{kn:deligne} Deligne, P. \emph{Cat\'egories tannakiennes}, The Grothendieck Festschrift,
  Vol.\ II, Progr. Math., vol.~87, Birkh\"auser Boston, Boston, MA, 1990,
  pp.~111--195.
\bibitem{kn:dt} Doi, Y. and  Takeuchi, M. 
{\em BiFrobenius algebras.} In Andruskiewtisch, N. (ed.) et al. 
{\em New trends in Hopf
  algebra theory.} Proceedings of the colloquium on quantum groups and
  Hopf algebras, La Falda, Sierras de C\'ordoba, Argentina, August 9-13,
  1999. Providence, RI: American Mathematical Society
  (AMS) (2000). 
\bibitem{kn:dri} Drinfel'd,V. G. {\em Quasi--Hopf algebras.}\, Leningrad
  Math. {\bf 1}, (1990), pp. 1419--1457.
\bibitem{kn:eno} Etingof, P., Nikshych, D. and Ostrik, V. {\em An analogue
  of Radofrd's $\mathcal S^4$ formula for finite tensor
  categories.}\,Int. Math. Res. Not. {\bf 54}, (2004), pp. 2915-2933.

\bibitem{kn:fs} Ferrer Santos, W. and Haim, M. {\em Radford's formula
    for biFrobenius algebras and applications.} To appear in
  Communications in Algebra. 
\bibitem{kn:hn} Hausser, F. and Nill, F. \emph{Integral Theory for Quasi-Hopf
  Algebras.} arXiv:math.QA/9904164.
\bibitem{kn:johnstone} Johnstone, P.T. \newblock {\em Sketches of an elephant: a topos theory compendium. {V}ol. 1},
  volume~43 of {\em Oxford Logic Guides}.
\newblock The Clarendon Press Oxford University Press, New York, 2002.

\bibitem{kn:JS} Joyal, A. and Street, R. \emph{Braided tensor
    categories.} Adv. Math, {\bf 102}, 1, (1993), pp. 20--78.
\bibitem{kn:ks} Kadison, L. and Stolin, A. A. {\em An approach to Hopf
  algebras via Frobenius coordinates I.} \, Beitr. Algebra
  Geom. {\bf 42}, 2,  (2001), pp. 359-384. 

\bibitem{kn:k} Kadison, L. {\em An approach to Hopf algebras via
  Frobenius Coordinates.}\,J. Algebra, {\bf 295}, 1, (2006), pp.
  27-43. 
\bibitem{kn:kassel} Kassel, C. {\em Quantum Groups.}\, New York:
  Springer--Verlag. Graduate Texts in Mathematics, {\bf 155}, 1995,
  xii+523 p. 
\bibitem{kn:kelly}
Kelly, G.M. \emph{Basic concepts of enriched category theory},
 London Mathematical Society Lecture Note Series, vol.~64, Cambridge
 University Press, Cambridge, 1982. 
\bibitem{kn:l} Larson, R. {\em Characters of Hopf algebras.} J. Algebra,
  {\bf 17}, (1971), pp. 352-368.
\bibitem{kn:nacho} Lopez Franco, I. {\em Hopf modules for autonomous
    pseudomonoids and the monoidal centre}.  (preprint)
  {\tt http://uk.arxiv.org/abs/0710.3853}.
\bibitem{kn:n} Nikshych, D. {\em On the structure of weak Hopf algebras.}
  Adv. Math, {\bf 170}, (2002), pp. 257-286.
\bibitem{kn:r} Radford, D. {\em The order of the antipode of a finite
  dimensional Hopf algebra.} Amer. J. Math. {\bf 98}, (1976), pp.
333-355.
\bibitem{kn:schau1} Schauenburg, P. \emph{Hopf modules and the double of a quasi-{H}opf
  algebra}, Trans. Amer. Math. Soc. \textbf{354} (2002), no.~8, 3349--3378
  (electronic).
\bibitem{kn:Schauenburg} Schauenburg, P. \emph{Two
    characterizations of finite quasi-{H}opf algebras.} J. Algebra,
  {\bf 273}, 2, (2004), pp. 538--550 (2004).
\bibitem{kn:schneider} Schneider, H.-J. {\em Lectures on Hopf
  Algebras}, notes by S. Natale. Trabajos de Matem\'atica, Vol. 31/95,
  (1995), FaMAF, C\'ordoba. Argentina.
\bibitem{kn:s} Sweedler, M. {\em Integrals for Hopf algebras.} Ann. of
  Math. {\bf 91},(1969), pp.  323-335. 
\bibitem{kn:sweedler} Sweedler, M. {\em Hopf algebras.} New York:
  W.A. Benjamin, Inc. 1969, 336 p. 
\end{thebibliography}
\end{document}